\documentclass[smallextended,envcountsect]{svjour3} 

\pdfoutput=1

\smartqed 
\usepackage{graphicx}
\usepackage{amsmath,amssymb}
\usepackage{amsfonts}
\usepackage{algorithm}
\usepackage{graphics}
\usepackage{epsfig}
\usepackage{enumitem}
\usepackage{algpseudocode}
\usepackage[dvipsnames]{xcolor}
\usepackage{subcaption}
\usepackage{wrapfig}
\usepackage{xcolor}
\usepackage{lipsum}
\usepackage{float}
\usepackage{comment}
\usepackage[mathscr]{eucal}
\usepackage{biblatex}
\newcommand{\stochgrad}[4]{\ensuremath{\widehat{\nabla}_{#1} #2(#3,#4)}}

\newcommand{\id}{\mathrm{Id}}

\newcommand{\reals}{\mathbb{R}}

\newcommand{\Exp}{\operatorname{\mathbb{E}}}
\newcommand{\fix}{\mathrm{fix}}
\newcommand{\zer}{\mathrm{zer}}

\newcommand{\basicmethod}{{basic method}}
\newcommand{\accalg}{{acceleration algorithm}}
\newcommand{\inpr}[2]{\ensuremath{\mathord{\left\langle #1, #2 \right\rangle}}}
\newcommand{\nat}{\ensuremath{\mathbb{N}}}
\newcommand{\prox}{\operatorname{prox}}
\newcommand{\grad}{\ensuremath{\mathrm{\nabla}}}
\newcommand{\z}{z}

\newcommand{\zaa}{\z}

\newcommand{\lv}{\left\Vert}
\newcommand{\rv}{\right\Vert}
\newcommand{\lp}{\left(}
\newcommand{\rp}{\right)}

\newcommand{\T}{\mathcal{T}}
\newcommand{\R}{{\mathcal{R}}}

\newcommand{\Rbar}{\Bar{\mathcal{R}}}

\newcommand{\norm}[1]{\ensuremath{\mathord{\left\Vert #1 \right\Vert}}}
\newcommand{\gra}[1]{\operatorname{gra}(#1)}
\newcommand{\seq}[1]{\ensuremath{(#1)_{k\in\mathbb{N}}}}

\newtheorem{assumption}{Assumption}[section]

\journalname{}

\begin{document}
\title{Hybrid Acceleration Scheme for Variance Reduced Stochastic Optimization Algorithms}
\author{Hamed Sadeghi   \and  Pontus Giselsson}
\institute{Hamed Sadeghi,  Corresponding author \at
             Department of Automatic Control, Lund University \\
             Lund, Sweden\\
             hamed.sadeghi@control.lth.se  
           \and
           Pontus Giselsson \at
             Department of Automatic Control, Lund University\\
             Lund, Sweden\\
             pontus.giselsson@control.lth.se
}
\date{Received: date / Accepted: date}
\maketitle

\begin{abstract}
Stochastic variance reduced optimization methods are known to be globally convergent while they suffer from slow local convergence, especially when moderate or high accuracy is needed. To alleviate this problem, we propose an optimization algorithm---which we refer to as a hybrid acceleration scheme---for a class of proximal variance reduced stochastic optimization algorithms. The proposed optimization scheme combines a fast locally convergent algorithm, such as a quasi--Newton method, with a globally convergent variance reduced stochastic algorithm, for instance SAGA or L--SVRG. Our global convergence result of the hybrid acceleration method is based on specific safeguard conditions that need to be satisfied for a step of the locally fast convergent method to be accepted.

We prove that the sequence of the iterates generated by the hybrid acceleration scheme converges almost surely to a solution of the underlying optimization problem. We also provide numerical experiments that show significantly improved convergence of the hybrid acceleration scheme compared to the basic stochastic variance reduced optimization algorithm.
\end{abstract}



\section{Introduction}\label{sec:intro}
We consider convex finite--sum optimization problems of the form
\begin{align} \label{eq:optProblem}
    \underset{x \in \reals^d}{\text{minimize}} \hspace{2mm} F(x)+g(x),
\end{align}
where $F:\reals^d\to\reals$ is the average of convex and smooth functions $f_i:\reals^d\to\reals$, {that is,}
\begin{align*}
    F(x) = \tfrac{1}{N} \sum_{i=1}^N f_i(x),
\end{align*}
for all $x\in\reals^d$ and $g:\reals^d\to\reals\cup\{\infty\}$ is  a  closed, convex,  proper, and potentially non--smooth function that can be used as a regularization term or to model convex constraints. Such finite--sum optimization problems are common in machine learning and statistics where they are known as regularized empirical risk minimization problems \cite{schmidt2017minimizing,teo2007scalable}. 

One approach to solve the finite--sum optimization problem~\eqref{eq:optProblem} is to use the proximal--gradient method \cite{beck2009fast}. However, at each iteration, the proximal--gradient algorithm requires as many individual gradient evaluations as the number of component functions of the finite--sum, which can be computationally expensive. Another approach is to apply stochastic proximal--gradient descent \cite{Nitanda2014StochasticPG,Rosasco2014ConvergenceOS}, which requires only one gradient evaluation at each iteration, but, due to the variance in the estimation of the full gradient, suffers from sub--linear convergence rate, even in the strongly convex setting \cite{johnson2013accelerating,kovalev2019don}. Several stochastic variance--reduced optimization algorithms such as SDCA \cite{shalev2013stochastic}, SVRG \cite{johnson2013accelerating}, and SAGA \cite{defazio2014saga}, have been designed to reduce the gradient approximation variance. These methods have been shown to be practically efficient and achieve global (linear) convergence for (strongly) convex problems. However, their local convergence is often slow in practice. 

To improve convergence, pre--determined data preconditioning \cite{li2017preconditioned,yang2016data} or metric selection \cite{giselsson2015metric} can be used. These are generic approaches that can be applied on top of acceleration schemes. However, finding the optimal or even a good metric is problem-- and algorithm--dependent and might be computationally expensive. Quasi--Newton type methods, such as Anderson acceleration \cite{Anderson1965IterativePF,walker2011anderson} and limited--memory BFGS \cite{liu1989limited}, instead find a suitable metric on the fly. Compared to stochastic optimization algorithms, these methods have higher per--iteration cost, but, often exhibit very fast local convergence. However, global convergence results are  scarce for non--smooth problems, whereas some results exist for fully smooth problems \cite{Rodomanov2021GreedyQM,Rodomanov2021NewRO,Rodomanov2021RatesOS}.

In this paper, we provide a generic algorithm that combines a method with locally fast convergence (that will be called {\emph{acceleration method}}) with a globally convergent proximal stochastic optimization algorithm (that will be called {\emph{basic method}}). The key feature of the general algorithm is a set of safeguard conditions that decide if an acceleration step can be accepted while maintaining global convergence. If the safeguard conditions are not satisfied, a step of the basic method is taken. This results in a hybrid scheme that automatically selects between two different algorithms and benefits both from the global convergence properties of the basic method and the fast local convergence of the acceleration method. We refer to our proposed algorithm as the \emph{hybrid acceleration scheme}.

The idea of a hybrid algorithm that selects between a globally convergent method and locally fast, but not globally, convergent method has been explored, e.g., in \cite{themelis2019supermann,zhang2018globally}, whose selection criteria are extensions of the one in \cite{giselsson2016line}. A key difference between our approach and \cite{themelis2019supermann,zhang2018globally} is that their methods are based on a deterministic basic method, while ours is based on a variance reduced stochastic method. This difference necessitates a completely different convergence analysis and enables for faster progress far from the solution in our finite--sum problem setting since our method takes advantage of that particular problem structure \cite{schmidt2017minimizing}. 

Due to the flexibility of our scheme, many different locally fast methods can be used. For instance; limited--memory BFGS (lBFGS) \cite{liu1989limited}, Anderson acceleration \cite{Anderson1965IterativePF}, and the class of \emph{vector extrapolation methods} \cite{smith1987extrapolation} to which, e.g., the regularized nonlinear acceleration  \cite{scieur2016regularized} and its stochastic counterpart \cite{scieur2017nonlinear} belong.

We instantiate our hybrid method with two different local methods, namely limited--memory BFGS (lBFGS) \cite{liu1989limited}  and Anderson acceleration \cite{Anderson1965IterativePF}. In our numerical experiments, we combine these methods with Loop--less SVRG \cite{kovalev2019don} in our hybrid acceleration method. Our numerical experiments show that our hybrid acceleration scheme can exhibit significantly improved convergence compared to the basic stochastic optimization algorithm.

The paper is outlined as follows. In Section~\ref{sec:Prelim}, we recall some basic definitions. Section~\ref{sec:formulation} discusses the problem formulation and the link between deterministic and stochastic gradient methods and introduces the family of stochastic optimization algorithms that is considered in this work. In Section~\ref{sec:propsed-alg}, the hybrid acceleration method is introduced and in Section~\ref{sec:analysis}, we prove its convergence. Numerical experiments are presented in Section~\ref{sec:numex} and concluding remarks are given in Section~\ref{sec:conclusion}.

\section{Preliminaries} \label{sec:Prelim}
The set of the real numbers and the $d$-dimensional Euclidean space are denoted by $\reals$ and $\reals^d$  respectively. For a symmetric positive definite matrix $\Gamma$ and $x, y \in \reals^d$, $\langle x,y\rangle$, $\| x \|$, and $\| x \|_\Gamma$ are the inner product, the induced norm, and the weighted norm $\|x\|_{\Gamma}:=\sqrt{\langle x,\Gamma x\rangle}$ respectively. Moreover, the $d \times d$ identity matrix is denoted by $I_d$.

The notation $2^{\reals^d}$ denotes the power set of $\reals^d$. A map $A:\reals^d\rightrightarrows 2^{\reals^d}$ is characterized by its graph $\gra{A} = \{(x,u)\in\reals^d\times\reals^d : u\in Ax\}$. The operator $A$ is monotone, if $\inpr{u-v}{x-y}\geq0$ for all $(x,u),(y,v)\in\gra{A}$. A monotone operator $A$ is maximally monotone if there exists no monotone operator $B:\reals^d\rightrightarrows  2^{\reals^d}$ such that $\gra{B}$ properly contains $\gra{A}$. A mapping $T:\reals^d \rightarrow \reals^d$ is $L$\emph{-Lipschitz continuous} if $\|T(x) - T(y) \| \leq L\|x - y \|$ for all $x, y \in \reals^d$, and is nonexpansive if it is 1-Lipschitz continuous. Further $T:\reals^d \rightarrow \reals^d$ is
\begin{enumerate}[label=\roman*)]
    \item \emph{firmly nonexpansive} if
        \begin{equation*}
             \|x - y - (T(x)- T(y)) \|^2 \leq \|x - y \|^2 - \|T(x) - T(y) \|^2  \quad\quad \forall x, y \in \reals^d,
        \end{equation*}
    \item $\tfrac{1}{L}$\emph{-cocoercive} if
        \begin{equation*}
            \langle T(x) - T(y) , x-y \rangle \geq \tfrac{1}{L}\|T(x) - T(y) \|^2  \quad\quad \forall x, y \in \reals^d.
        \end{equation*}
\end{enumerate}
For a mapping $T$, $\tfrac{1}{L}$-cocoercivity implies its $L$-Lipschitz continuity. The other direction does not hold in general. However, if the mapping is the gradient of a convex function, then its $L$-Lipschitz continuity and $\tfrac{1}{L}$-cocoercivity are equivalent {\cite[Corollary 18.17]{bauschke2017convex}}. A differentiable function $F: \reals^d \rightarrow \reals$ is said to be $L$-\emph{smooth}, if its gradient is $L$-Lipschitz continuous.

The subdifferential of a function $f:\reals^d\rightarrow\reals\cup\{\infty\}$ at $x\in\reals^d$ is denoted by $\partial f(x)$ and defined as
\begin{align*}
    \partial f(x) = \{v\in\reals^d : f(y)\geq f(x) +\inpr{v}{y-x} \text{~for all~} y\in\reals^d \}.
\end{align*}
The proximal mapping of a closed, convex and proper function $g: \reals^d \rightarrow \reals \cup \{+\infty \}$, is defined as $\prox_{\lambda g} (v) = \underset{x}{\mathrm{argmin}} \left(g(x) + \tfrac{1}{2\lambda}\| x-v \|^2 \right)$, where $\lambda>0$.

The set of fixed--points of a mapping $\T:\reals^d \rightarrow \reals^d$, is denoted by $\fix\T$ and defined as $\fix\T = \{ x\in\reals^d : x = \T x \}$. The zero--set of a map $\R:\reals^d \rightarrow \reals^d$ is indicated by $\zer\R$ and given by $\zer\R = \{ x\in\reals^d : 0=\R x \}$. It is evident that $\fix \T = \zer{(\id-\T)}$, where $\id-\T$ is the residual map of the operator $\T$.

\section{Problem formulation and basic method}\label{sec:formulation}
We are interested in solving the following convex optimization problem
\begin{align} \label{eq:optProblem-repeat}
    \underset{x \in \reals^d}{\text{minimize}} \hspace{2mm} \tfrac{1}{N} \sum_{i=1}^N f_i(x)+g(x),
\end{align}
under the following assumptions.
\begin{assumption} \label{ass:assumptions1}
We assume that
\renewcommand{\labelenumi}{\emph{(\roman{enumi})}}
\begin{enumerate}
    \item For {each} $i\in \{1,\ldots,N \}$, the function $f_i: \reals^d \rightarrow \reals$ is convex, differentiable and $L_i$-smooth.
    \item The function $g: \reals^d \rightarrow \reals \cup \{+\infty \}$ is convex, closed and proper.
    \item The solution set of the problem is nonempty.
\end{enumerate}
\end{assumption}

The necessary and sufficient optimality condition for this problem is given by Fermat's rule as
\begin{equation} \label{eq:inclusion1}
    0\in{\partial (F+g)(x)=}\grad F(x)+\partial g(x){,}
\end{equation}
{where the equality holds since all $f_i$ have full domain and $g$ is proper {\cite[Theorem 16.3 and Corollary 16.48]{bauschke2017convex}}}. {This} means that any $x^\star$ that satisfies the optimality condition~{\eqref{eq:inclusion1}}, is a solution to the associated optimization problem~{\eqref{eq:optProblem-repeat}}. It is also known that fixed--points of the {proximal--gradient} operator, namely, the set $\{x\in\reals^d : x=\prox_{\lambda g}(x - \lambda \nabla F(x)), \lambda>0\}$, are solutions of  problem~\eqref{eq:optProblem-repeat}. In fact, all solutions of the inclusion problem~\eqref{eq:inclusion1}, are fixed--points of the proximal--gradient operator or, equivalently, zeros of its residual mapping, which is given by 
\begin{equation*}
        \R x = x -  \prox_{\lambda g}\left(x - \lambda \nabla F(x) \right),
\end{equation*}
for any $\lambda > 0$ {\cite[Section 4.2]{parikh2014proximal}}. For $0<\lambda<\tfrac{2}{L}$ with $L$ being the smoothness modulus of $F$, iterating the proximal gradient mapping finds a solution of problem~\eqref{eq:optProblem-repeat} \cite{combettes2011proximal}.

The optimality condition~\eqref{eq:inclusion1}, can be reformulated in a primal--dual form by storing all gradients of component functions{ $f_i$}. In that case, the optimality condition {becomes}
\begin{align} \label{eq:inclusion2}
\begin{cases}
    &0\in \partial g(x)+\tfrac{1}{N}\sum_{i=1}^N y_i \\
    &0= y_1 - \grad f_1(x)\\
    &\quad\vdots\\
    &0= y_N - \grad f_N(x)
\end{cases},
\end{align}
where $y_i$ denotes the $i$--th dual variable. This is clearly equivalent to \eqref{eq:inclusion1}. Therefore, a primal--dual solution $\z^\star:=(x^\star,y_1^\star,\ldots,y_N^\star)$ satisfies \eqref{eq:inclusion2}, if and only if $x^\star$ satisfies \eqref{eq:inclusion1} and is a solution of \eqref{eq:optProblem-repeat}. It also holds that $\z^\star$ satisfies \eqref{eq:inclusion2} if and only if it satisfies $\Rbar\z^\star=0$, where $\Rbar$ is the primal--dual residual mapping
\begin{equation} \label{eq:primal--dual-res-map}
    \Rbar\z :=
            \begin{pmatrix}
                x - \prox_{\lambda g} \lp x - \tfrac{\lambda}{N}\sum_{i=1}^N y_i \rp \\
                   y_1 - \grad f_1(x)  \\
                \vdots \\
                   y_N - \grad f_N(x) 
            \end{pmatrix},
\end{equation}
in which $\z = (x,y_1,\ldots,y_N)$ is the primal--dual variable. We record the equivalence between zeroes of $\Rbar$ and solutions to \eqref{eq:optProblem-repeat}, in Proposition~\ref{prop:p_pd_equivalence} (with proof in Appendix~\ref{app:prop3.1}) and Lipschitz continuity of $\Rbar$ in Proposition~\ref{prop:Lipschitz-cont} (with proof in Appendix~\ref{app:prop3.2}).

\begin{proposition}\label{prop:p_pd_equivalence}
Given the residual map in \eqref{eq:primal--dual-res-map}, the primal--dual point $\z^\star=(x^\star,y_1^\star,\ldots,y_N^\star)$ satisfies $\Rbar z^\star=0$, if and only if $x^\star$ solves \eqref{eq:optProblem-repeat}. Furthermore, for each index $i$, $y_i^\star$ is unique. 
\end{proposition} 

\begin{proposition} \label{prop:Lipschitz-cont}
Let for all $i \in \{1,\ldots,N\}$, $ f_i(x)$ be $L_i$-smooth, then the primal--dual residual mapping, $\Rbar$ in \eqref{eq:primal--dual-res-map}, is Lipschitz-continuous.
\end{proposition}

In order to find zeros of $\Rbar$, one way is to form iterates based on \eqref{eq:primal--dual-res-map}, and  evaluate all $y_i$'s (the full gradient) at each iteration, which would be {similar to the} proximal--gradient algorithm. However, when $N$ is very large, a key challenge is the high per--iteration cost of $N$ gradient evaluation{s} which makes the algorithm {very} expensive. This gives rise to the idea of using a cheaply evaluable approximation of the true gradient instead, and randomly evaluate gradients of only one or some of $f_i$'s at each iteration. The following gives such an approximation
\begin{equation} \label{eq:stochApprox}
     \stochgrad{i_k}{F}{x}{y} \triangleq \tfrac{1}{N p_{i_k}} \left( \grad f_{i_k}(x) - y_{i_k}\right) + \tfrac{1}{N} \sum\limits_{i=1}^N y_i ,
\end{equation}
in which $i_k$ is an index randomly drawn from $\{1,\ldots,N \}$ based on some probability distribution and $p_{i_k}$ is its associated probability. This stochastic approximation is based on the average of the dual variables which is modified by a correction term, $( \grad f_{i_k}(x) - y_{i_k})/({N p_{i_k}})$. The correction term is added in order to progressively improve the approximation by incorporating the latest gradient information and also to make $\stochgrad{i_k}{F}{x}{y}$ an unbiased estimate of the true gradient. 

Using the approximation $\grad F(x) \approx \stochgrad{i_k}{F}{x}{y}$, and inspired by the proximal gradient algorithm, a family of proximal stochastic optimization algorithms can be formulated as
\begin{equation} \label{eq:basicmethod1}
\begin{aligned} 
    x^{k+1} &= \prox_{\lambda g} \left( x^k - \lambda \stochgrad{i_k}{F}{x^k}{y^k} \right),\\
    y_i^{k+1} &= y_i^k + \varepsilon_i^k (\grad f_i (x^k) - y_i^k),  \quad  \forall i \in \{ 1, \hdots, N \},
\end{aligned}
\end{equation}
where $k$ is the iteration counter, $x^k$ is the primal variable, $y^k = (y_1^k,\ldots,y_N^k)$ with $y_i^k$ being the $i$--th dual variable, $\lambda>0$ is the step size, $\varepsilon_i^k {\in \{ 0,1\}}$ is a random binary variable that determines whether the $i$--th dual variable is to be updated at iteration $k$ (the associated probability of $\varepsilon_i^k = 1$ is $\rho_i$), and $\stochgrad{i_k}{F}{x^k}{y^k}$ is {the} stochastic approximation of the true gradient that is defined in \eqref{eq:stochApprox}. This approximation of the full gradient is unbiased since 
\begin{align*}
    \Exp_k(\stochgrad{i_k}{F}{x^k}{y^k}) &= \tfrac{1}{N}\sum_{i=1}^N(\grad f_i(x^k) - y_i^k) + \tfrac{1}{N}\sum_{i=1}^N y_i^k\\
    &= \tfrac{1}{N}\sum_{i=1}^N\grad f_i(x^k) = \grad F(x^k){,}
\end{align*}
in which, $\Exp_k$ denotes expected value operation given all available information up to step $k$.  We refer to \eqref{eq:basicmethod1} as the \emph{basic method}. On the other hand, there are algorithms that use a biased estimation of the true gradient \cite{morin2019svag,roux2012stochastic}, but in this work we only consider the unbiased case. The algorithm in \eqref{eq:basicmethod1} has been analyzed in \cite{davis2016smart} in the monotone operator setting and in \cite{morin2020sampling} in the strongly convex setting.

The class of stochastic optimization algorithms~\eqref{eq:basicmethod1} has L--SVRG \cite{davis2016smart,kovalev2019don} and SAGA \cite{defazio2014saga}  as special cases. The L--SVRG algorithm is extracted from  \eqref{eq:basicmethod1} with uniform sampling of $i_k\in\{1,\ldots,N \}$ and
\begin{align*}
    \varepsilon_i^k = \begin{cases} 1  \qquad \mathrm{if} \hspace{1mm} q<\rho \\ 0 \qquad\mathrm{otherwise}\end{cases}, \quad \forall i\in \{1,\ldots,N\},
\end{align*}
 where $q$ is uniformly sampled from $[0,1]$ and $0<\rho\leq1$. Therefore, all dual variables are updated together and on average once every $\rho^{-1}$ iterations. The algorithm in \eqref{eq:basicmethod1} reduces to SAGA with $i_k$ uniformly sampled from $\{1,\hdots, N\}$ and
 \begin{align*}
    \varepsilon_i^k = \begin{cases} 1  \qquad \mathrm{if} \hspace{1mm} i_k = i \\ 0 \qquad\mathrm{otherwise}\end{cases} , \quad \forall i\in \{1,\ldots,N\}.
\end{align*}
Therefore, for SAGA, at each iteration, only one of the dual variables is updated and the others remain unchanged.

\section{Hybrid acceleration scheme} \label{sec:propsed-alg}
In this section, we introduce a novel hybrid strategy to accelerate local convergence of proximal stochastic optimization algorithms of the form \eqref{eq:basicmethod1}, in which {the} approximation of {the} true gradient and the update law of the dual variables, vary depending on the choice of basic method. The basic method~\eqref{eq:basicmethod1}, is devised to solve large--scale finite--sum optimization problems of the form \eqref{eq:optProblem-repeat}, and is globally convergent while it has slow local convergence. Therefore, in our acceleration scheme, they are combined with a locally fast convergent method. The proposed acceleration scheme is given in Algorithm~\ref{alg:proposed} and discussed below.

\begin{algorithm}
\caption{\hspace{1mm} General framework of the hybrid acceleration scheme}
	\begin{algorithmic}[1]
		\State \textbf{Input:} initial point $\z^0$, positive constants $C$, $D$, $\delta$, merit function $V(.)$, acceleration algorithm $\mathcal{A}(.)$ and its memory size $m$ (if needed), $K_0$, the basic method and its parameters, primal and dual probability distribution, the step size $\lambda$, $\Gamma=\mathrm{blkdiag}(I_d,\tfrac{\lambda}{N\rho_1 L_1} I_d,\hdots,\tfrac{\lambda}{N\rho_N L_N} I_d)$, and the maximum permissible number of iterations $k_{max}$.
		\State set $k = k_{\mathrm{aa}}=0$
		\While{$k<k_{max}$}
		\State $m_{k_{\mathrm{aa}}} = \min \{m,k_{\mathrm{aa}}\}$
		\State find $\zaa^+=\mathcal{A}(\z^k,\zaa^{{k_{\mathrm{aa}}-1}},\ldots,\zaa^{{k_{\mathrm{aa}}-m_{k_{\mathrm{aa}}}}})$ from acceleration algorithm
        \If{$V( \zaa^+ )    \leq \frac{C V(\z^0)}{(k_{\mathrm{aa}}+1)^{(1+\delta)}}$ \textbf{and} $\lv \zaa^+ - \zaa^k \rv_\Gamma \leq D V(\zaa^{k})$}
        \State set $\z^{k+1} = \zaa^+$ and $k_{\mathrm{aa}} \leftarrow k_{\mathrm{aa}}+1$
        \State $k\leftarrow k+1$ and $k_{\mathrm{aa}} \leftarrow k_{\mathrm{aa}}+1$
        \Else
        \State set $(\widetilde{x}^{0},\widetilde{y}_1^{0},\ldots,\widetilde{y}_N^{0}) = z^k$
        \For{$s=0, 1, \cdots, K_0-1$}
        {\begin{equation*} \begin{aligned}
		&\begin{cases}
		    \widetilde{x}^{s+1} = \prox_{\lambda g}  \left( \widetilde{x}^s - \lambda \stochgrad{i_s}{F}{\widetilde{x}^s}{\widetilde{y}^s}\right), \\
            \widetilde{y}_i^{s+1} = \widetilde{y}_i^s + \varepsilon_i^s (\grad f_i (\widetilde{x}^s) - \widetilde{y}_i^s), \quad \forall i \in \{1 ,\hdots, N  \}.
        \end{cases}\\
        &\hspace{2.8mm} \z^{k+1} = (\widetilde{x}^{s+1},\widetilde{y}_1^{s+1},\ldots,\widetilde{y}_N^{s+1})\\
        &\hspace{2.8mm}k\leftarrow k+1
        \end{aligned}\end{equation*}}
        \EndFor
        \EndIf
		\EndWhile
	\end{algorithmic}
\label{alg:proposed}
\end{algorithm}

\paragraph{\textbf{Description of the algorithm.}}
In order to initialize the scheme one needs to select (i) the parameters and probability distributions used in the basic method; (ii) an acceleration algorithm $\mathcal{A}(.)$ along with its associated parameters; and (iii) an initial point $\z^0$. The acceleration algorithm $\mathcal{A}(.)$ can be algorithms such as lBFGS or Anderson acceleration that both store and use a history of past $m$ iterates to find a next iterate. Then, the algorithm works as follows: at the beginning of each iteration the iterate from the acceleration algorithm, $\zaa^+$, has to be computed. If $\zaa^+$ satisfies some \textit{safeguard conditions}, that we will discuss below, we set it as the true next iterate, $\z^{k+1}$, and the main counter of the loop, $k$, and also the acceleration algorithm's counter, $k_{aa}$, are increased by one; then, we proceed to the next iteration. Otherwise, $K_0$ steps of the basic method are performed in the inner loop of the algorithm. It is evident that $K_0$ can differ among different iterations of the outer loop of the scheme, but, we considered it as a constant for simplicity. The algorithm is to be run as above until the last iteration is reached or some termination criteria are met. Note that if the iterate from the acceleration algorithm is accepted, the basic method steps need not to be performed, that is, the basic method and the acceleration algorithm are not being run in parallel.

\paragraph{\textbf{Safeguard conditions and merit function.}}
For a nominal next iterate of {the} acceleration algorithm, $\zaa^+$, to be accepted as the actual next iterate of the scheme, the following conditions have to be satisfied
\begin{align}
    V( \zaa^+ )    &\leq C V(\z^0)(1+k_{\mathrm{aa}})^{-(1+\delta)}  \label{eq:safeguard1}, \\
    \lv \zaa^+ - \z^k \rv_\Gamma &\leq D V(\zaa^{k})  \label{eq:safeguard2},
\end{align}
where $\delta$, $C$ and $D$ are positive constants, $V:\reals^{(N+1)d} \rightarrow \reals$ is a merit function (that is discussed below), and
\begin{equation*}
    \Gamma=\mathrm{blkdiag}(I_d,\tfrac{\lambda}{N\rho_1 L_1} I_d,\hdots,\tfrac{\lambda}{N\rho_N L_N} I_d).
\end{equation*}
The safeguard condition~\eqref{eq:safeguard1} enforces the merit function  to be convergent to zero. Condition~\eqref{eq:safeguard2}, is to ensure that the sequence $(\|\zaa^{k+1} - \zaa^k\|_\Gamma)_{k\in I_{\mathrm{aa}}}$, where $I_{\mathrm{aa}}$ is the set of indices for which the next iterate is obtained from \text{\accalg}, is diminishing and finally convergent to zero.

Our convergence theory, that will be given in the next section, assumes that; i) the merit function outputs nonnegative values, ii) for any sequence $\seq{z^k}$, the merit function is such that
\begin{align*}
    V(\z^k)\to 0 \qquad\Longrightarrow\qquad \lv\Rbar \z^k\rv\to 0.
\end{align*}
Therefore, a feasible choice for the merit function can be the following scaled $l_2$--norm of $\Rbar \z^{k}$
\begin{equation}\label{eq:merit_func}
    V(\z^k) = \lv \Rbar\z^k \rv_\Gamma.
\end{equation}
For this choice of merit function, which we use in this work, both requirements on the merit function are met. Other options for the merit function could be the sum or maximum of the vector of {the} last $p$ scaled $l_2$--norm of residuals.

\section{Convergence results}\label{sec:analysis}
In this section, we provide results on convergence of the basic method and the hybrid acceleration scheme. {B}efore proceeding to convergence results, we summarize the notations and the assumptions that are used in the theorems and their proof{s}. The proofs are given in the Appendix.
  
\paragraph{Notation.} $\mathcal{X}^\star$ indicates the solution set of problem~\eqref{eq:optProblem-repeat}, $\z = (x,y_1,\ldots,y_N)$ denotes a primal--dual variable for \eqref{eq:primal--dual-res-map}, $\Rbar$ is the primal--dual residual operator defined in \eqref{eq:primal--dual-res-map}, $\z^\star = (x^\star,y_1^\star,\ldots,y_N^\star)$ is an arbitrary point in the set of zeros of {the} primal--dual residual mapping with $y_i^\star = \grad f_i (x^\star)$,  $p_i$ is primal sampling probability, $\rho_i$ is the $i$--th dual variable update probability,  $\lambda>0$ is the step size, $\Exp_k$ denotes the expected value operator given all the information up to the $k$--th iteration, and  $\Gamma=\mathrm{blkdiag}(I_d,\tfrac{\lambda}{N\rho_1 L_1} I_d,\hdots,\tfrac{\lambda}{N\rho_N L_N} I_d)$. Moreover, $\z^k = (x^k,y_1^k,\ldots,y_N^k)$ denotes {the} $k$--th primal--dual iterate; and  $\seq{z^k}$, $\seq{x^k}$, and $\seq{y_i^k}$ are the sequences of primal--dual-, primal-, and the $i$--th dual iterates, respectively. 

The following is a result that is used in proof of Theorem~\ref{thm:main}. The proof can be found in Appendix~\ref{app:prop5.1}.

\begin{proposition} \label{prop:dual-convergence}
Under Assumption~\ref{ass:assumptions1}, almost sure (a.s.) convergence of $\seq{\z^k}$ to a $\Bar{\z} \in \zer{\Rbar}$, implies a.s. convergence of $\seq{x^k}$ and $\seq{y_i^k}$ to a $\Bar{x} \in \mathcal{X}^\star$ and $\Bar{y}_i = \grad f_i(\Bar{x})$ respectively. 
\end{proposition}

The result in Theorem~\ref{thm:basicMethod} and its proof (given in Appendix~\ref{app:theorem5.1}) share similarities with \cite{davis2016smart, morin2020sampling}.

\begin{theorem}\label{thm:basicMethod}
Let $\z^k$ be the $k$--th primal--dual iterate associated with the basic method iterates in \eqref{eq:basicmethod1}, then given Assumption~\ref{ass:assumptions1}, the following holds
\begin{equation*}
    \Exp_k \Vert \z^{k+1} - \z^\star\Vert_\Gamma^2  \leq   \Vert \z^{k} - \z^\star\Vert_\Gamma^2 - \zeta_k,
\end{equation*}
with 
\begin{equation}\label{eq:zeta_seq}
    \begin{aligned}
    \zeta_k &= \sum_{i=1}^N \tfrac{\lambda}{N}(\tfrac{1}{L_i} - \tfrac{2\lambda}{N p_i}) \left(\Vert \grad f_i(x^k) - \grad f_i(x^\star) \Vert^2 +  \Vert y_i^k - y_i^\star \Vert^2 \right)\\
    &\qquad + \tfrac{2\lambda^2}{N^2} \Vert\sum_{i=1}^N(y_i^k-y_i^\star)\Vert^2+\lambda^2\Vert\grad F(x^k)-\grad F(x^\star)\Vert^2\\
    &\qquad + \Exp_k \Vert x^{k+1}-x^k+\lambda (\stochgrad{i_k}{F}{x^k}{y^k} - \grad F(x^\star)) \Vert^2.
\end{aligned}
\end{equation}
 Furthermore, if $0<\lambda<\mathrm{min}_i\{\tfrac{N p_i}{2 L_i}\}$, then $\seq{\z^k}$ converges a.s. to  a random variable $\Bar{z} \in \zer{\Rbar}$ and $\seq{x^k}$ and $\seq{y_i^k}$ converge a.s. to random variables $\Bar{x} \in 
\mathcal{X}^\star$ and $\Bar{y}_i = \grad f_i (\Bar{x})$  respectively.
\end{theorem}

\begin{remark}
In order to ensure a.s. convergence in Theorem~\ref{thm:basicMethod}, the coefficient of all terms in $\zeta_k$ must be positive. Then, from relation~\eqref{eq:zeta_seq}, it is evident that for each $i$, $0<\lambda < \tfrac{N p_i}{2 L_i}$ must hold. Therefore, the {smallest of these} has to be set as the upper bound of {$\lambda$}, that is, $0<\lambda<\mathrm{min}_i\{\tfrac{N p_i}{2 L_i}\}$. The largest {upper bound} of the step size is attained when we have Lipschitz probability distribution for primal sampling, namely, $p_i = \tfrac{L_i}{\sum_{i=1}^N L_i}$.
\end{remark}

The following result is on a.s. convergence of the sequence of iterates that are obtained from Algorithm~\ref{alg:proposed}. The proof is presented in Appendix~\ref{app:theorem5.2}.

\begin{theorem}\label{thm:main}
Suppose that Assumption~\ref{ass:assumptions1} holds, that $0<\lambda<\mathrm{min}_i\{\tfrac{N p_i}{2 L_i}\}$, and that {the merit function $V:\reals^{(N+1)d}\to\reals$ is nonnegative and such that for all sequences $\seq{z^k}$ satisfying $V(\z^k)\to 0$ we have $\|\Rbar\z^k\|\to 0$}. Then $\seq{\z^k}$ in Algorithm~\ref{alg:proposed} converges a.s. to a random variable $\Bar{z} \in \zer{\Rbar}$. Moreover, $\seq{x^k}$ and $\seq{y_i^k}$ converge a.s. to random variables  $\Bar{x} \in 
\mathcal{X}^\star$ and $\Bar{y}_i = \grad f_i (\Bar{x})$  respectively.
\end{theorem}

\section{Numerical experiments}\label{sec:numex}
We solve a regularized logistic regression problem for binary classification of the form
\begin{align}\label{eq:logReg}
    \underset{x = (w,b)}{\text{minimize}} \hspace{2mm} \sum_{i=1}^N {\text{log}(1 + e^{\theta_i^T w + b}) - u_i (\theta_i^T w + b)} + \tfrac{\xi}{2}\Vert w \Vert_2^2,
\end{align}
where $\theta_i \in \reals^d$ and $u_i \in \{ 0,1 \}$ are training data and labels respectively, and $\xi > 0$ is a regularization parameter. The optimization problem variable is $x = (w,b)$ with $w \in \reals^d$ and $b \in \reals$.

In the hybrid acceleration scheme, we use L--SVRG as the basic method and either Anderson acceleration or lBFGS as the acceleration algorithms. The following lists the algorithms that are used in the numerical experiments
\begin{itemize}
    \item[$\bullet$] GD: Gradient descent method with fixed step size,
    \item[$\bullet$] L--SVRG: Loopless Stochastic Variance Reduced Gradient method,
    \item[$\bullet$] L--SVRG+AA: L--SVRG as the \text{\basicmethod} combined with Anderson acceleration,
    \item[$\bullet$] L--SVRG+lBFGS: L--SVRG as the \text{\basicmethod} combined with limited--memory BFGS.
\end{itemize}

In order to use Anderson acceleration as the acceleration algorithm in the hybrid acceleration scheme, an associated fixed--point mapping of problem~\eqref{eq:logReg} is needed. Let $F(x)$ denote the objective function of problem~\eqref{eq:logReg}. Then, the associated mapping of the problem that is used by Anderson acceleration is given by
\begin{align*}
    \T_{gd}(x) = x - \lambda\nabla F(x)
\end{align*}
for $\lambda\in(0,2/L)$, where $L$ is smoothness modulus of $F$. On the other hand, since the objective function at hand  has no non--smooth part, the lBFGS algorithm can also be utilized in the hybrid acceleration scheme to solve problem~\eqref{eq:logReg}. Unlike Anderson acceleration, lBFGS method does not need an associated fixed--point map of the problem, rather, it requires gradients of the objective function in order to find a solution. See Appendix~\ref{subsec:Anderson} and Appendix~\ref{subsec:lBFGS} for descriptions of Anderson acceleration and lBFGS methods, respectively.

A rough approximation of the per--iteration count of floating point operations for the different algorithms are as follows
\begin{itemize}
    \item[$\bullet$] $4Nd$ for gradient descent,
    \item[$\bullet$] $12Nd$ for L--SVRG,
    \item[$\bullet$] $4Nd+\tfrac{4}{3}m^3+2m^2d$ for Anderson acceleration,
    \item[$\bullet$] $4Nd+2d^2+13md + \xi_{bt}4Nd$ for lBFGS,
\end{itemize}
where, $N$ is the number of the component functions (which is the same as the number of samples in the training dataset), $d$ is the dimension of the optimization problem variable, $m$ is the size of memory stack for either Anderson acceleration or lBFGS and $\xi_{bt}$ is a coefficient to include an approximate average cost for back-tracking line search of lBFGS.

Numerical simulations are done using two datasets; \emph{UCI Madelon}  \cite{CC01a} with 2000 samples and 500 features, and \emph{UCI Sonar} \cite{CC01a} with 208 samples and 60 features. In the numerical experiments, the regularization parameter in the objective function is set to $\xi=0.01$, we used a memory size of $m=5$ for both Anderson acceleration and lBFGS, the constants of the safeguard condition of the hybrid acceleration scheme is $C=D=10^6$, and $\delta = 10^{-6}$, and the parameter $K_0$ is set to the number of samples of the associated dataset. Moreover, we used the merit function as is defined in \eqref{eq:merit_func}.

\begin{figure}
\begin{subfigure}{0.49\textwidth}
\centering
\includegraphics[ width = 0.99\linewidth ,keepaspectratio]{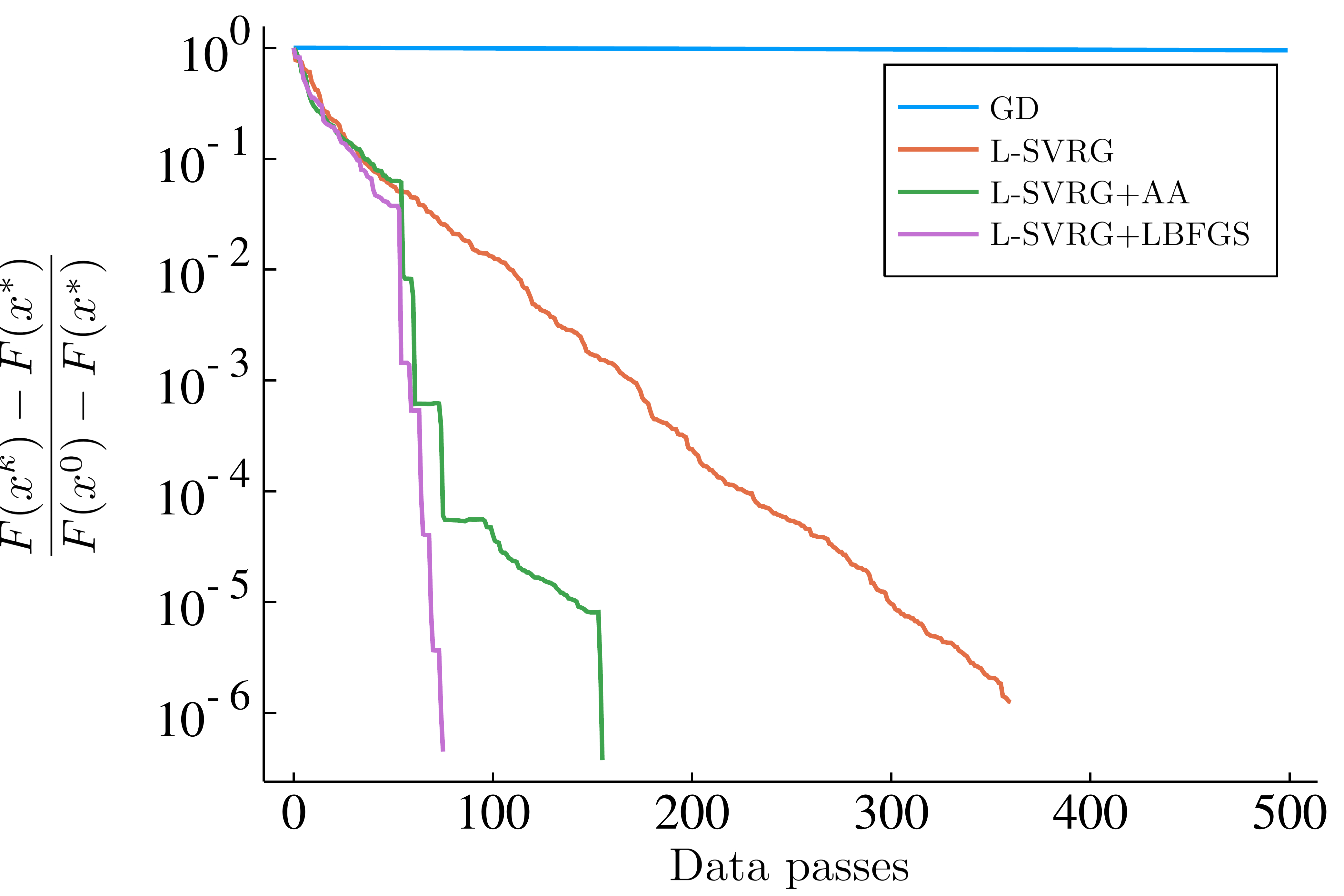}
\end{subfigure}
\begin{subfigure}{0.49\textwidth}
\centering
\includegraphics[ width = 0.99\linewidth ,keepaspectratio]{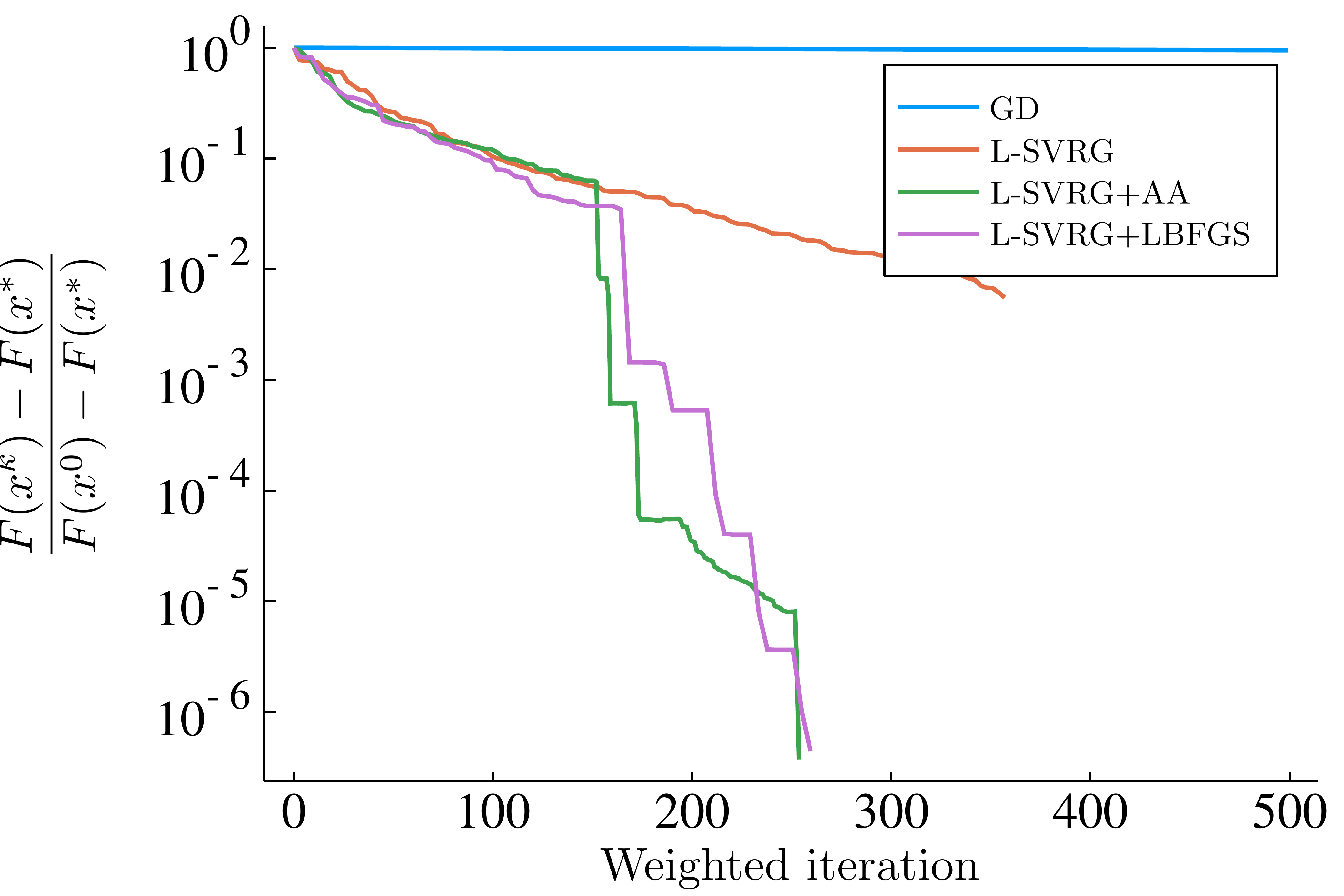}
\end{subfigure}   
\caption{Normalized sub-optimality vs. number of passes over data (the plot to the left) and weighted iteration number (the plot to the right) for the logistic regression problem \eqref{eq:logReg}, on \emph{UCI Madelon} dataset (2000 samples, 500 features), solved using GD, L--SVRG, L--SVRG+AA and L--SVRG+lBFGS methods with regularization parameter $\xi = 0.01$.}
\label{fig:madelon}
\end{figure}

\begin{figure}
\begin{subfigure}{0.49\textwidth}
\centering
\includegraphics[ width = 0.99\linewidth ,keepaspectratio]{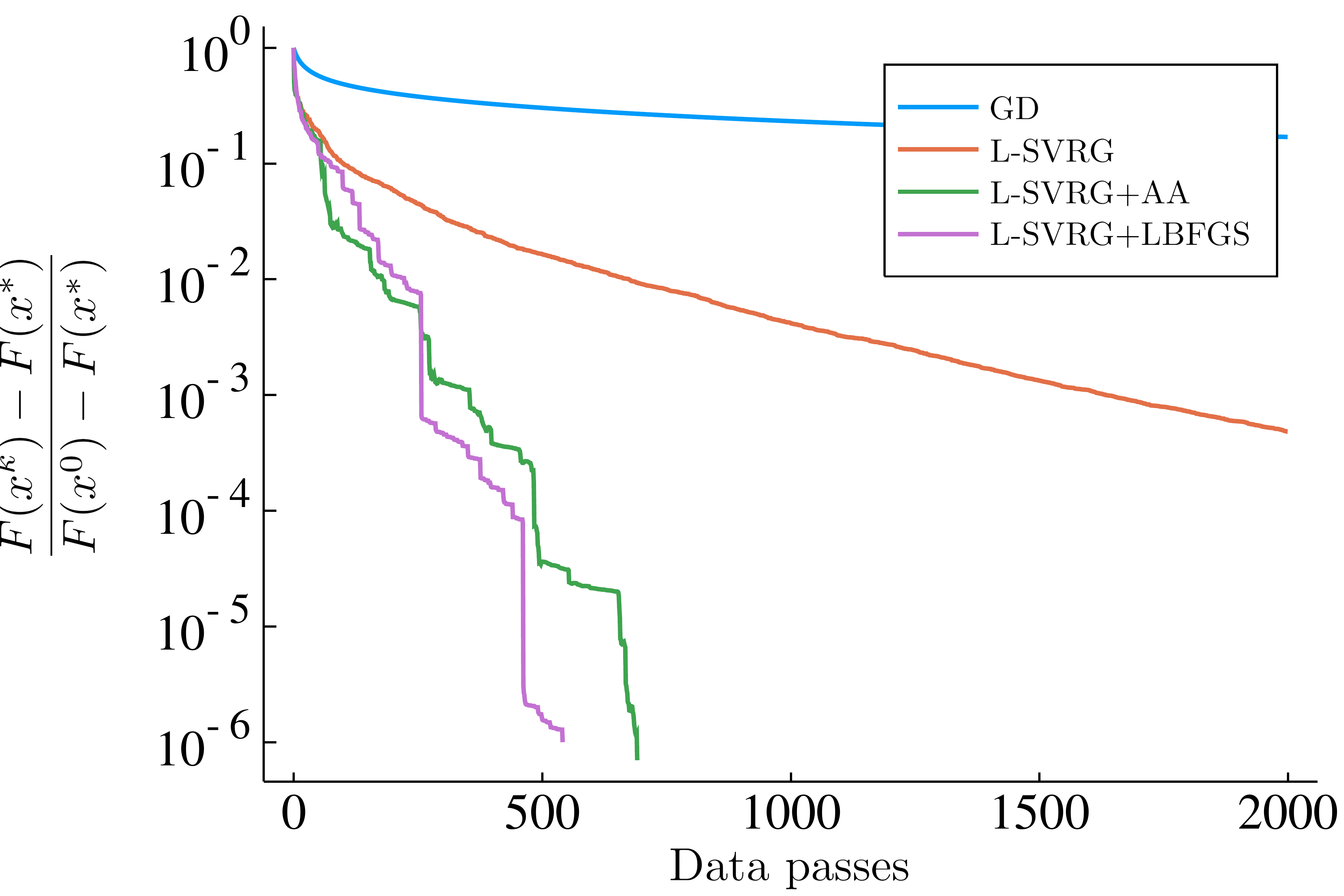}
\end{subfigure}
\begin{subfigure}{0.49\textwidth}
\centering
\includegraphics[ width = 0.99\linewidth ,keepaspectratio]{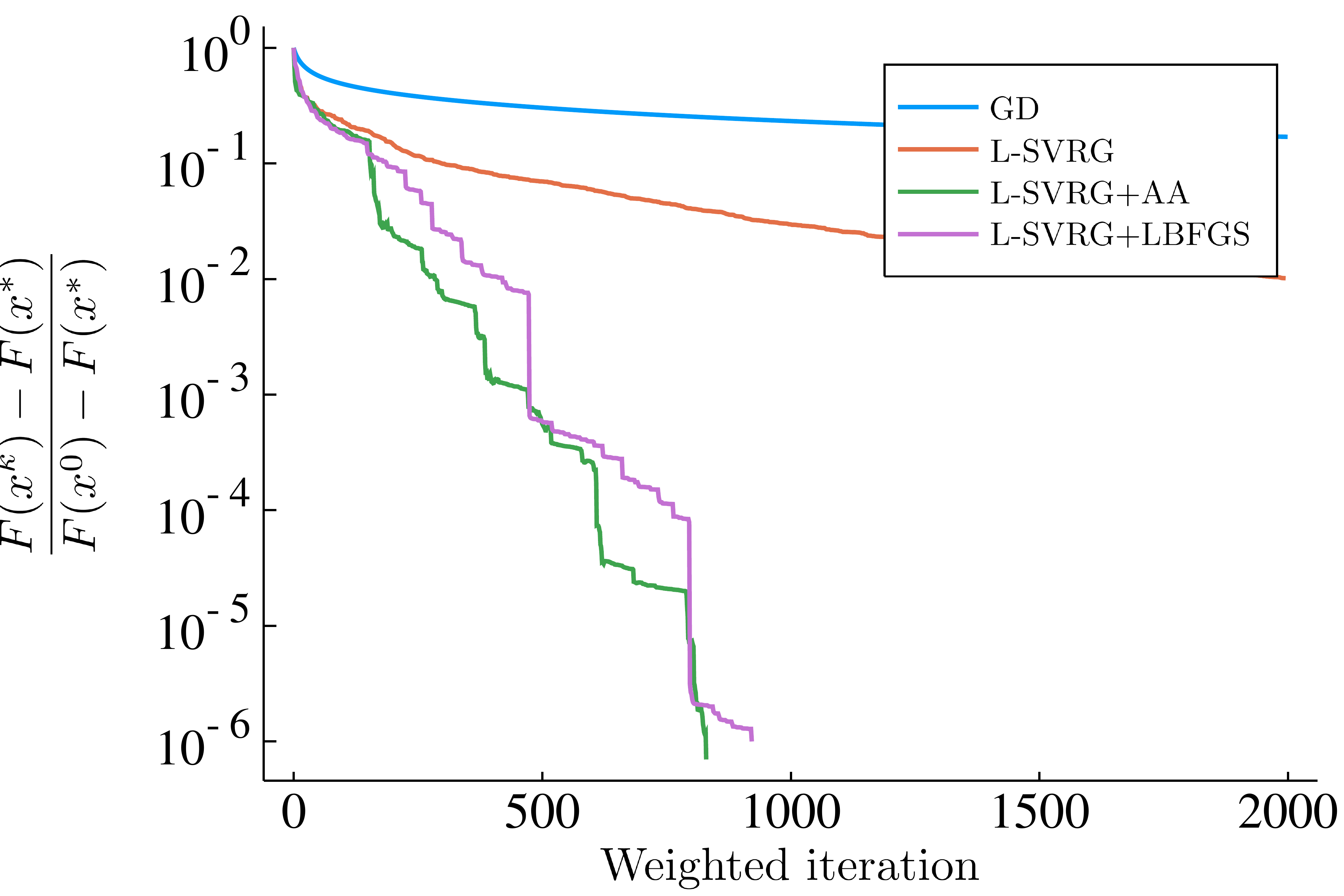}
\end{subfigure}
\caption{Normalized sub-optimality vs. number of passes over data (the plot to the left) and weighted iteration number (the plot to the right) for the logistic regression problem \eqref{eq:logReg}, on \emph{UCI Sonar} dataset (208 samples, 60 features), solved using GD, L--SVRG, L--SVRG+AA and L--SVRG+lBFGS methods with regularization parameter $\xi = 0.01$.}
\label{fig:sonar}
\end{figure}

In Figure~\ref{fig:madelon} and Figure~\ref{fig:sonar}, the left plots show relative value of objective function versus step number (which is basically equal to the  total number of full gradient evaluations up to that step), and the right plot illustrates the relative value of objective function versus weighted iteration counts. The weighted iteration is intended to include a rough approximation of computational cost in such a way that different methods at each weighted iteration have roughly the same computational expense. Therefore, it provides a better comparison in terms of computational complexity among different algorithms. The simulation results show remarkable improvement in convergence rate and overall computational cost of the hybrid acceleration scheme compared to those of the basic method.

\section{Conclusion}\label{sec:conclusion}
In this paper, we proposed and showed almost sure convergence of a hybrid acceleration scheme. It combines a globally convergent variance reduced stochastic gradient method---the basic method---with a fast locally convergent method---the acceleration method---to benefit from the strengths of both methods; global convergence of the basic method and fast local convergence of the acceleration method. Our numerical experiments show that our algorithm performs significantly better than the basic method in isolation, while preserving global convergence guarantees that the local acceleration methods lack.

\begin{acknowledgements}
This work was partially supported by the Wallenberg AI, Autonomous Systems and Software Program (WASP) funded by the Knut and Alice Wallenberg Foundation.

The authors would like to thank Bo Bernhardsson for his fruitful feedback on this work.
\end{acknowledgements}

\appendix
\section*{Appendix A}
\renewcommand{\thesubsection}{A.\arabic{subsection}}
\renewcommand{\thelemma}{A.\arabic{lemma}}

In what follows, we provide the proofs of the propositions and the theorems that are not addressed in the body of the paper.

\subsection{Proof of Proposition~\ref{prop:p_pd_equivalence}} \label{app:prop3.1}

From $\Rbar \z^\star = 0$ and for any $\lambda>0$ we have
\begin{align*}
      \begin{cases} x^\star - \prox_{\lambda g} \lp x^\star - \tfrac{\lambda}{N}\sum_{i=1}^N y_i^\star \rp=0\\ y_1^\star = \grad f_1(x^\star)\\\qquad\vdots \\
        y_N^\star = \grad f_N(x^\star) \end{cases} &\quad\iff\quad x^\star -  \prox_{\lambda g}\left(x^\star - \lambda \nabla F(x^\star) \right)=0\\
        &\quad\iff\quad 0 \in \nabla F(x^\star) + \partial g(x^\star)\\
        &\quad\iff\quad 0 \in \partial (F+g)(x^\star),
\end{align*}
where the last equivalence holds due to $f_i$'s and $g$ having full domain \cite[Theorem 16.3 and Corollary 16.48]{bauschke2017convex}. Therefore, by Fermat's rule $x^\star$ is a solution of problem~\eqref{eq:optProblem-repeat}. 

Now suppose that ${x}_1^\star$ and ${x}_2^\star$ are two distinct solutions to the problem, that is
\begin{align*}
    -\nabla F({x}_1^\star)&\in\partial g({x}_1^\star), \\
    -\nabla F({x}_2^\star)&\in\partial g({x}_2^\star).
\end{align*}
Then using the fact that $\partial g$ is monotone and that each $\nabla f_i$ is $\tfrac{1}{L_i}$-cocoercive, we have
\begin{align*}
     0\geq\left\langle {x}_2^\star - {x}_1^\star,    \grad F \lp {x}_2^\star \rp -  \grad F \lp {x}_1^\star \rp \right\rangle \geq \sum_{i=1}^N \tfrac{1}{N L_i}\lv \grad f_i ({x}_2^\star) - \grad f_i({x}_1^\star) \rv^2 \geq 0,
\end{align*}
which gives that ${y}_i^\star = \grad f_i({x}_2^\star) = \grad f_i({x}_1^\star)$ for all $i$'s. Hence it follows that ${y}_i^\star = \grad f_i ({x}^\star)$ is unique.
\qed

\subsection{Proof of Proposition~\ref{prop:Lipschitz-cont}} \label{app:prop3.2}

In the following proof, we use nonexpansiveness of the proximal operator and $L_i$-Lipschitz continuity of $\grad f_i(x)$ for all $i$
\begingroup
\allowdisplaybreaks
\begin{align*}
    \norm{\Rbar\hat{\z}-\Rbar\z}^2 &=  \norm{\hat{x} - \prox_{\lambda g} \bigg(\hat{x} - \tfrac{\lambda}{N}\sum_{i=1}^N{\hat{y}_i}\bigg) - x + \prox_{\lambda g} \bigg(x-\tfrac{\lambda}{N}\sum_{i=1}^N{y_i} \bigg)}^2 \\
    &\qquad+ \sum_{i=1}^N{\norm{\hat{y}_i - \grad f_i(\hat{x}) - y_i + \grad f_i (x)}}^2\\
    &\leq \Bigg(\norm{\hat{x}-x} + \norm{\hat{x} - \tfrac{\lambda}{N}\sum_{i=1}^N{\hat{y}_i} - \lp{x - \tfrac{\lambda}{N}\sum_{i=1}^N{{y_i}}}\rp}\Bigg)^2 \\
    &\qquad + \sum_{i=1}^N 2\lp\norm{\hat{y}_i-y_i}^2 + \norm{\grad f_i(\hat{x}) - \grad f_i(x)}^2\rp\\
    &\leq 2\norm{\hat{x}-x}^2 + 2\norm{\hat{x} - \tfrac{\lambda}{N}\sum_{i=1}^N{\hat{y}_i} - \lp{x - \tfrac{\lambda}{N}\sum_{i=1}^N{{y_i}}}\rp}^2 \\
    &\qquad + \sum_{i=1}^N 2\lp\norm{\hat{y}_i-y_i}^2 + \norm{\grad f_i(\hat{x}) - \grad f_i(x)}^2\rp\\
    &\leq  6\norm{\hat{x}-x}^2  + 4\tfrac{\lambda^2}{N^2}\norm{\sum_{i=1}^N{(\hat{y}_i-y_i)}}^2\\
    &\qquad + \sum_{i=1}^N 2\lp\norm{{\hat{y}_i-y_i}\rv^2 + L_i^2\lv{\hat{x} - x}}^2\rp\\
    &\leq {(6 + 2\sum_{i=1}^N{L_i^2})}\norm{\hat{x}-x}^2 + {(2 + 4\tfrac{\lambda^2}{N})}\sum_{i=1}^N{\norm{\hat{y}_i-y_i}^2}\\
    &\leq \Bar{\alpha} \lp{\norm{\hat{x}-x}^2 + \sum_{i=1}^N{\norm{\hat{y}_i-y_i}^2} }\rp = \Bar{\alpha} \norm{\hat{\z}-\z}^2,
\end{align*}
\endgroup
where
\begin{align*}
    \Bar{\alpha} = \max \left( 6 + 2\sum_{i=1}^N{L_i^2},  2 + 4\tfrac{\lambda^2}{N}\right).
\end{align*}
The first inequality is given by the triangle inequality and nonexpansiveness of proximal operator for the first term and by Young's inequality for terms in the sum. The second and third inequalities is given by Young's inequality, and the fourth one is given by $\|a_1+\ldots+a_N\|_2^2 \leq N(\|a_1\|_2^2+\ldots+\|a_N\|_2^2)$. Therefore, $\Rbar$ is Lipschitz continuous.
\qed

\subsection{Proof of Proposition~\ref{prop:dual-convergence}} \label{app:prop5.1}

Using the definition of {the} primal--dual residual operator at $\Bar{\z} = (\Bar{x},\Bar{y}_1.\ldots,\Bar{y}_N)$ 
\begin{equation*}
    \Rbar \Bar{\z} \triangleq
            \begin{pmatrix}
                \Bar{x} - \prox_{\lambda g} \lp \Bar{x} - \tfrac{\lambda}{N}\sum_{i=1}^N \Bar{y}_i \rp \\
                   \Bar{y}_1 - \grad f_1(\Bar{x})  \\
                \vdots \\
                   \Bar{y}_N - \grad f_N(\Bar{x}) 
            \end{pmatrix} = 0
\end{equation*}
gives $\Bar{y}_i = \grad f_i (\Bar{x})$ for all $i$. This in turn yields
\begin{align*}
    \Bar{x} - \prox_{\lambda g} \lp\Bar{x} - \tfrac{\lambda}{N}\sum_{i=1}^N \grad f_i(\Bar{x})\rp = \Bar{x} - \prox_{\lambda g} \lp\Bar{x} - \lambda\grad F(\Bar{x})\rp = 0,
\end{align*}
which evidently means that $\Bar{x} \in \mathcal{X}^\star${, i.e., $0\in\partial g(\Bar{x})+\nabla F(\Bar{x})$}. Since $z^k$ converges to $\bar z$ almost surely, $x^k$ and $y_i^k$, respectively, converge to $\bar x$ and $\Bar{y}_i=\nabla f_i(\bar x)$ almost surely. This concludes the proof.
\qed

\subsection{Lemmas for proof of Theorem~\ref{thm:basicMethod}} \label{app:lemmas-theorem5.1}

The following lemmas are needed in our proof of Theorem~\ref{thm:basicMethod}.

\begin{lemma} {\label{lemma1}}
Let $\Rbar$ be the primal--dual residual operator \eqref{eq:primal--dual-res-map}, $(x^\star,y^\star)\in \zer{\Rbar}$, $f_i$ be convex and $L_i$-Lipschitz continuous, and $y_i^\star = \grad f_i(x^\star)$ for all $i \in \{1, ..., N  \}$ , then for the iterates given in \eqref{eq:basicmethod1}, the following bounds {the} variance of the primal variable
\begin{equation} \label{eq:lemma1}
    \begin{aligned}
        \Exp_k \Vert x^{k+1} - x^\star \Vert^2 &\leq \Vert x^k - x^\star \Vert^2  - \sum_{i=1}^N {\tfrac{2\lambda}{NL_i} \Vert \grad f_i(x^k) - \grad f_i(x^\star) \Vert^2}\\
     &\qquad   - \Exp_k \Vert x^{k+1} - x^k + \lambda (\stochgrad{i_k}{F}{x^k}{y^k} - \grad F(x^\star)) \Vert^2\\
    &\hspace{33mm}+  \lambda^2 \Exp_k \Vert \stochgrad{i_k}{F}{x^k}{y^k} - \grad F(x^\star) \Vert^2.
    \end{aligned}
\end{equation}
\end{lemma}

{\it Proof.} Using firm nonexpansiveness of the proximal operator, we have
\begingroup
\allowdisplaybreaks
\begin{align*}
    \Vert x^{k+1} - x^\star \Vert^2 &= \left\Vert \prox_{\lambda g} \left(x^k - \lambda \stochgrad{i_k}{F}{x^k}{y^k}\right) - \prox_{\lambda g} \left(x^\star - \lambda \grad F(x^\star)\right) \right\Vert^2\\
        &\leq  \left\Vert  x^k - \lambda \stochgrad{i_k}{F}{x^k}{y^k} - (x^\star - \lambda \grad F(x^\star))\right\Vert^2 - \\
        & \quad\quad\quad\quad \Vert x^k - \lambda \stochgrad{i_k}{F}{x^k}{y^k} - \prox_{\lambda g}(x^k - \lambda \stochgrad{i_k}{F}{x^k}{y^k}) \\
        & \quad\quad\quad\quad\quad\quad - (x^\star - \lambda \grad F(x^\star)) + \prox_{\lambda g}(x^\star - \lambda \grad F(x^\star)) \Vert^2 \\
        &= \Vert x^k - x^\star\Vert^2 - 2\lambda\langle x^k-x^\star, \stochgrad{i_k}{F}{x^k}{y^k} - \grad F(x^\star)\rangle \\
        &\qquad + \lambda^2 \Vert \stochgrad{i_k}{F}{x^k}{y^k} - \grad F(x^\star) \Vert^2\\ 
        &\qquad - \Vert x^{k+1} - x^k + \lambda (\stochgrad{i_k}{F}{x^k}{y^k} - \grad F(x^\star)) \Vert^2.
\end{align*}
\endgroup
The second equality  above, is given by $x^\star=\prox_{\lambda g}(x^\star - \lambda \grad F(x^\star))$ and by the primal update formula $x^{k+1} = \prox_{\lambda g}(x^k - \lambda \stochgrad{i_k}{F}{x^k}{y^k})$. Taking expected value conditioned on all available information up to step $k$, yields
\begingroup
\allowdisplaybreaks
\begin{align*}
    \Exp_k \Vert x^{k+1} - x^\star \Vert^2 &\leq \Vert x^k - x^\star\Vert^2 - 2\lambda\langle x^k-x^\star, \grad{F}(x^k) - \grad F(x^\star)\rangle\\
    &\qquad\qquad - \Exp_k \Vert x^{k+1} - x^k + \lambda (\stochgrad{i_k}{F}{x^k}{y^k} - \grad F(x^\star)) \Vert^2\\
    &\qquad\qquad\qquad\qquad+ \lambda^2\Exp_k\Vert\stochgrad{i_k}{F}{x^k}{y^k} - \grad F(x^\star)\Vert^2\\
    &\leq \Vert x^k - x^\star \Vert^2 - 2\lambda \sum_{i=1}^N {\tfrac{1}{NL_i} \Vert \grad f_i(x^k) - \grad f_i(x^\star) \Vert^2}\\
    &\qquad\qquad - \Exp_k \Vert x^{k+1} - x^k + \lambda (\stochgrad{i_k}{F}{x^k}{y^k} - \grad F(x^\star)) \Vert^2\\
    &\qquad\qquad\qquad\qquad + \lambda^2 \Exp_k \Vert \stochgrad{i_k}{F}{x^k}{y^k} - \grad F(x^\star) \Vert^2.
\end{align*}
\endgroup
In the first inequality, we used $\Exp_k \stochgrad{i_k}{F}{x^k}{y^k} = \grad F(x^k)$ and  the second inequality  is given by cocoercivity of $\grad f_i(x)$.
\qed

\begin{lemma}{\label{lemm2}}
Let $\Rbar$ be the primal--dual residual operator \eqref{eq:primal--dual-res-map}, $(x^\star,y^\star)\in \zer{\Rbar}$ and $y_i^\star = \grad f_i(x^\star)$ for all $i \in \{1, ..., N  \}$, then for the iterates given in \eqref{eq:basicmethod1}, the following holds:
\begin{equation} \label{eq:lemma2}
    \begin{aligned}
        \Exp_k \left( \sum_{i=1}^N {\tfrac{\lambda}{N \rho_i L_i} \Vert y_i^{k+1} - y_i^\star \Vert^2} \right) &= \sum_{i=1}^N \tfrac{\lambda}{N L_i} \Vert \grad f_i(x^k) - \grad f_i(x^\star) \Vert^2\\
     &\qquad\qquad  + \sum_{i=1}^N (1-\rho_i) \tfrac{\lambda}{N \rho_i L_i} \Vert y_i^k - y_i^\star\Vert^2 ,
    \end{aligned}
\end{equation}
where $\rho_i$ is the probability of $\varepsilon_i^k$ being $1$ for $y_i$.
\end{lemma}

{\it Proof.}
By substitution of $y_i^{k+1}$ from \eqref{eq:basicmethod1} we get 
\begingroup
\allowdisplaybreaks
\begin{align*}
    \Exp_k \left( \sum_{i=1}^N {\tfrac{\lambda}{N \rho_i L_i} \Vert y_i^{k+1} - y_i^\star \Vert^2} \right) &= \Exp_k \left( \sum_{i=1}^N \tfrac{\lambda}{N \rho_i L_i} \Vert y_i^k + \varepsilon_i^k(\grad f_i(x^k) - y_i^k) - y_i^\star \Vert^2 \right) \\
    &= \sum_{i=1}^N \tfrac{\lambda}{N \rho_i L_i}\Exp_k \Vert y_i^k + \varepsilon_i^k(\grad f_i(x^k) - y_i^k) - y_i^\star \Vert^2 \\
    &= \sum_{i=1}^N \tfrac{\lambda}{N \rho_i L_i}\Big(\rho_i\Vert\grad f_i(x^k) - y_i^\star\Vert^2 + (1-\rho_i)\Vert y_i^k  - y_i^\star \Vert^2\Big)\\
    &= \sum_{i=1}^N \tfrac{\lambda}{N L_i} \Vert \grad f_i(x^k) - \grad f_i(x^\star) \Vert^2\\
    &\qquad\quad  + \sum_{i=1}^N (1-\rho_i) \tfrac{\lambda}{N \rho_i L_i} \Vert y_i^k - y_i^\star\Vert^2.
\end{align*}
\endgroup
In the third equality we used the fact that the only random variable in the expression to the right of the second equality is $\varepsilon_i^k\in\{0,1\}$ and the probability of $\varepsilon_i^k$ being $1$ is assumed to be $\rho_i$.
\qed

\begin{lemma}\label{lemma3}
Let $\Rbar$ be the primal--dual residual operator of the problem, $(x^\star, y^\star)\in \zer{\Rbar}$  and $y_i^\star = \grad f_i(x^\star)$ for all $i \in \{1, ..., N  \}$, then for the iterates given in \eqref{eq:basicmethod1}, the following gives the update variance bound:
\begin{equation} \label{eq:lemma3}
\begin{aligned}
   \Exp_k\Vert\stochgrad{i_k}{F}{x^k}{y^k} - \grad F(x^\star) \Vert^2 \leq& \sum_{i=1}^N {\tfrac{2}{N^2 p_i} \lp \Vert\grad f_i(x^k)- \grad f_i(x^\star)\Vert^2+ \Vert y_i^k - y_i^\star\Vert^2 \rp}\\
   &- 2 \Vert \tfrac{1}{N}\sum_{i=1}^N {(y_i^k - y_i^\star)}\Vert^2 -\Vert\grad F(x^k) - \grad F(x^\star)\Vert^2.
   \end{aligned}
\end{equation}
\end{lemma}

{\it Proof.} We start with the left-hand side of \eqref{eq:lemma3}. Using the identity $\Exp \|X \|^2 = \|\Exp X \|^2 + \Exp \| X - \Exp X\|^2$, gives
\begin{equation}\label{eq:var_bound}
    \begin{aligned} 
   \Exp_k\Vert\stochgrad{i_k}{F}{x^k}{y^k} - \grad F(x^\star) \Vert^2 &=  \Vert\grad F(x^k) - \grad F(x^\star)\Vert^2 \\
   &\qquad+ \Exp_k \Vert\stochgrad{i_k}{F}{x^k}{y^k} - \grad F(x^k)\Vert^2.
\end{aligned}
\end{equation}
Now for the second term in the right-hand side, substitution of $\stochgrad{i_k}{F}{x^k}{y^k}$ yields
\begingroup
\allowdisplaybreaks
\begin{align*}
   \Exp_k\Big\Vert\stochgrad{i_k}{F}{x^k}{y^k} &- \grad F(x^k)\Big\Vert^2\\
   &= \Exp_k \norm{\tfrac{1}{N p_{i_k}} (\grad f_{i_k}(x^k) - y_{i_k}^k) + \tfrac{1}{N}\sum_{i=1}^N y_i^k - \grad F(x^k)}^2 \\
   &= \Exp_k\Bigg\Vert\tfrac{1}{N p_{i_k}} (\grad f_{i_k}(x^k) - \grad f_{i_k}(x^\star) + y_{i_k}^\star -  y_{i_k}^k) \\
   &\qquad\qquad\qquad + \tfrac{1}{N}\sum_{i=1}^N y_i^k - \tfrac{1}{N}\sum_{i=1}^N y_i^\star + \grad F(x^\star) - \grad F(x^k)\Bigg\Vert^2 \\
   &\leq 2\Exp_k \norm{\tfrac{1}{N p_{i_k}} (\grad f_{i_k}(x^k) - \grad f_{i_k}(x^\star)) - (\grad F(x^k) - \grad F(x^\star))}^2 \\
   &\qquad\qquad+ 2\Exp_k\norm{\tfrac{1}{N p_{i_k}}(y_{i_k}^k - y_{i_k}^\star)   - (\tfrac{1}{N}\sum_{i=1}^N y_i^k - \tfrac{1}{N}\sum_{i=1}^N y_i^\star)}^2\\
   &= 2\Exp_k \norm{\tfrac{1}{N p_{i_k}} (\grad f_{i_k}(x^k) - \grad f_{i_k}(x^\star))}^2 - 2\norm{\grad F(x^k) - \grad F(x^\star)}^2\\
   &\qquad\qquad+ 2\Exp_k\norm{\tfrac{1}{N p_{i_k}}(y_{i_k}^k - y_{i_k}^\star)}^2   - 2\norm{\tfrac{1}{N}\sum_{i=1}^N (y_i^k - y_i^\star)}^2\\ 
   &= 2\sum_{i=1}^N\tfrac{1}{N^2 p_{i}} \norm{\grad f_{i}(x^k) - \grad f_{i}(x^\star)}^2 - 2\norm{\grad F(x^k) - \grad F(x^\star)}^2\\
   &\qquad\qquad+ 2\sum_{i=1}^N\tfrac{1}{N^2 p_{i}}\norm{y_{i}^k - y_{i}^\star}^2  - 2\norm{\tfrac{1}{N}\sum_{i=1}^N (y_i^k - y_i^\star)}^2.
   \end{align*}
\endgroup
The inequality above is given by Cauchy-Schwarz and Young's inequalities. The third equality is given by the identity $\Exp \| X - \Exp X\|^2=\Exp \|X \|^2 - \|\Exp X \|^2$. Substituting in \eqref{eq:var_bound} yields
\allowdisplaybreaks
\begin{equation*}
    \begin{aligned}
       &\Exp_k\norm{\stochgrad{i_k}{F}{x^k}{y^k}-\grad F(x^\star)}^2\\
        &\hspace{10mm}\leq \sum_{i=1}^{N}\tfrac{2}{N^2 p_{i}} \norm{\grad f_{i}(x^k)-\grad f_{i}(x^\star)}^2 - \norm{\grad F(x^k)-\grad F(x^\star)}^2\\
        &\hspace{10mm}\qquad\qquad + \sum_{i=1}^{N}\tfrac{2}{N^2p_i}\norm{y_i^k-y_{i}^\star}^2-2\norm{\tfrac{1}{N}\sum_{i=1}^{N}(y_i^k-y_i^\star)}^2.
    \end{aligned}
\end{equation*}
\qed

\subsection{Proof of Theorem~\ref{thm:basicMethod}}\label{app:theorem5.1}
We first use the definition of $\Gamma$:
\begin{align}\label{eq:z-scaledNorm}
    \Vert z^k-z^\star\Vert_\Gamma^2 = \Vert x^k - x^\star\Vert^2 + \sum_{i=1}^N{\tfrac{\lambda}{N\rho_i L_i} \Vert y_i^k -y_i^\star\Vert^2}.
\end{align}
Then, adding \eqref{eq:lemma1} to \eqref{eq:lemma2} and reordering the terms, yield
\begin{align*}
    &\Exp_k\Vert x^{k+1} - x^\star\Vert^2 + \Exp_k\left(\sum_{i=1}^N{\tfrac{\lambda}{N\rho_i L_i} \Vert y_i^{k+1} -y_i^\star\Vert^2} \right) \\
    &\hspace{10mm}\leq \Vert x^k - x^\star\Vert^2   - \sum_{i=1}^N {\tfrac{2\lambda}{NL_i} \Vert \grad f_i(x^k) - \grad f_i(x^\star) \Vert^2} \\
    &\hspace{10mm}\qquad + \lambda^2 \Exp_k \Vert \stochgrad{i_k}{F}{x^k}{y^k}- \grad F(x^\star) \Vert^2\\
    &\hspace{10mm}\qquad  - \Exp_k \Vert x^{k+1} - x^k + \lambda (\stochgrad{i_k}{F}{x^k}{y^k} - \grad F(x^\star)) \Vert^2\\   
    &\hspace{10mm}\qquad+ \sum_{i=1}^N \tfrac{\lambda}{N L_i} \Vert \grad f_i(x^k) - \grad f_i(x^\star) \Vert^2  + \sum_{i=1}^N (1-\rho_i) \tfrac{\lambda}{N \rho_i L_i} \Vert y_i^k - y_i^\star\Vert^2\\
    &\hspace{10mm}= \Vert x^k - x^\star\Vert^2 + \sum_{i=1}^N  \tfrac{\lambda}{N \rho_i L_i} \Vert y_i^k - y_i^\star\Vert^2 - \sum_{i=1}^N {\tfrac{2\lambda}{NL_i} \Vert \grad f_i(x^k) - \grad f_i(x^\star) \Vert^2} \\
    &\hspace{10mm}\qquad + \lambda^2 \Exp_k \Vert \stochgrad{i_k}{F}{x^k}{y^k}- \grad F(x^\star) \Vert^2\\
    &\hspace{10mm}\qquad  - \Exp_k \Vert x^{k+1} - x^k + \lambda (\stochgrad{i_k}{F}{x^k}{y^k} - \grad F(x^\star)) \Vert^2\\   
    &\hspace{10mm}\qquad+ \sum_{i=1}^N \tfrac{\lambda}{N L_i} \Vert \grad f_i(x^k) - \grad f_i(x^\star) \Vert^2  - \sum_{i=1}^N \tfrac{\lambda}{N L_i} \Vert y_i^k - y_i^\star\Vert^2\\
\end{align*}
Now, we use \eqref{eq:z-scaledNorm} and \eqref{eq:lemma3} in the above inequality, which gives
\begingroup
\allowdisplaybreaks
\begin{align*}
    \Exp_k\Vert\z^{k+1} - \z^\star\Vert_\Gamma^2  &\leq \Vert\z^{k} - \z^\star\Vert_\Gamma^2 -  \sum_{i=1}^N {\tfrac{2\lambda}{NL_i} \Vert \grad f_i(x^k) - \grad f_i(x^\star) \Vert^2} - \lambda^2\Vert\grad F(x^k) - \grad F(x^\star)\Vert^2\\
    &\quad\quad + \sum_{i=1}^N \tfrac{\lambda}{N  L_i} \Vert \grad f_i(x^k) - \grad f_i(x^\star) \Vert^2 + \sum_{i=1}^N {(\tfrac{2\lambda^2}{N^2 p_i} - \tfrac{\lambda}{N  L_i})\Vert y_i^k - y_i^\star \Vert^2}\\
    &\quad\quad + \sum_{i=1}^N {\tfrac{2 \lambda^2}{N^2 p_i}\Vert \grad f_i(x^k) - \grad f_i(x^\star)\Vert^2} - \tfrac{2\lambda^2}{N^2} \Vert \sum_{i=1}^N {(y_i^k - y_i^\star)} \Vert^2\\
    &\quad\quad - \Exp_k \Vert x^{k+1} - x^k + \lambda (\stochgrad{i_k}{F}{x^k}{y^k} - \grad F(x^\star)) \Vert^2  \\ 
    &= \Vert\z^{k} - \z^\star\Vert_\Gamma^2 -  \sum_{i=1}^N {( \tfrac{\lambda}{NL_i} - \tfrac{2\lambda^2}{N^2 p_i} ) \Vert \grad f_i(x^k) - \grad f_i(x^\star) \Vert^2}\\
    &\quad\quad - \sum_{i=1}^N {(\tfrac{\lambda}{N L_i} - \tfrac{2\lambda^2}{N^2 p_i})\Vert y_i^k - y_i^\star \Vert^2} - \tfrac{2\lambda^2}{N^2} \Vert\sum_{i=1}^N {(y_i^k - y_i^\star)} \Vert^2 \\
    &\quad\quad - \Exp_k \Vert x^{k+1} - x^k + \lambda (\stochgrad{i_k}{F}{x^k}{y^k} - \grad F(x^\star)) \Vert^2 \\
    &\quad\quad- \lambda^2\Vert\grad F(x^k) - \grad F(x^\star)\Vert^2 \\
    &= \Vert\z^{k} - \z^\star\Vert_\Gamma^2 - \zeta_k
\end{align*}
\endgroup
where
\begingroup
\allowdisplaybreaks
\begin{align*}
    \zeta_k &= \sum_{i=1}^N {(\tfrac{\lambda}{NL_i} - \tfrac{2\lambda^2}{N^2 p_i}) \Vert\grad f_i(x^k) - \grad f_i(x^\star) \Vert^2}\\
    &\quad\quad + \sum_{i=1}^N {(\tfrac{\lambda}{N L_i} - \tfrac{2\lambda^2}{N^2 p_i})\Vert y_i^k - y_i^\star \Vert^2} + \tfrac{2\lambda^2}{N^2} \Vert\sum_{i=1}^N {(y_i^k - y_i^\star)} \Vert^2 \\
    &\quad\quad + \Exp_k \Vert x^{k+1} - x^k + \lambda (\stochgrad{i_k}{F}{x^k}{y^k} - \grad F(x^\star)) \Vert^2 + \lambda^2\Vert\grad F(x^k) - \grad F(x^\star)\Vert^2.
\end{align*}
\endgroup
This proves the first part of the theorem. To show a.s. convergence of $\seq{x^k}$ to a random variable in $\mathcal{X}^\star$, in view of {\cite[Proposition 2.3]{combettes2015stochastic}}, we need to show that the set of  sequential cluster points of the sequence $\seq{x^k}$ is a subset of $\mathcal{X}^\star$, then a.s. convergence of $\seq{x^k}$ to an $\mathcal{X}^\star$-valued random variable will follow. In the following, all limits and convergences are to be considered to hold almost surely, also if it is not explicitly written.

We choose $\lambda$ such that $0<\lambda<\mathrm{min}_i\{\tfrac{N p_i}{2 L_i}\}$ holds. This choice of $\lambda$, enforces non-negativeness to all the coefficients in relation \eqref{eq:zeta_seq}; thus, we have  $\zeta_k\geq0$ for all $k\in\nat$. Now using {\cite[Proposition 2.3.i]{combettes2015stochastic}}, we get that $\seq{\zeta_k}$ is a.s. summable. It follows by a.s. summability of $\seq{\zeta_k}$ that  both $\seq{y_i^k}$ and $\seq{\grad f_i(x^k)}$ converge to $\grad f_i(x^\star)$ almost surely. This in turn means that, as $k\to\infty$, $\stochgrad{i_k}{F}{x^k}{y^k} \rightarrow \grad F(x^\star)$ almsot surely.  Moreover, a.s. summability of $\seq{\zeta_k}$ implies that $\seq{ \Exp_k (\Vert x^{k+1}-x^k+\lambda(\stochgrad{i_k}{F}{x^k}{y^k}-\grad F(x^\star))\Vert^2) }$ a.s. converges to zero as $k \rightarrow \infty$ and since $\stochgrad{i_k}{F}{x^k}{y^k} - \grad F(x^\star) \rightarrow 0$,  we have that $\Exp_k (\Vert x^{k+1} - x^k\Vert^2 ) \rightarrow 0$, which implies $x^{k+1} - x^k \rightarrow 0$ almost surely.
Now, since the Euclidean space $\reals^{(N+1)d}$ is separable and $\zer{\Rbar}$ is closed, using {\cite[Proposition 2.3.iii]{combettes2015stochastic}}, for every $\z^\star \in \zer{\Rbar}$, the sequence $\seq{\Vert \z^k - \z^\star \Vert}$ converges almost surely. Summability of $\seq{\zeta_k}$ implies that $\Vert\z^{k} - \z^\star\Vert_\Gamma^2 - \Vert x^k - x^\star \Vert^2 \rightarrow 0$, and therefore, we infer that for every $x^\star \in \mathcal{X}
^\star$, the sequence $\seq{\Vert x^k - x^\star \Vert^2}$ is a.s. convergent, and therefore, the sequence $\seq{x^k}$ is bounded. Boundedness of $\seq{x^k}$ implies that it has at least one convergent subsequence. Denote this subsequence by $\seq{x^{n_k}}$. Now, from the optimality condition of the proximal operator we get
\begin{align*}
    0 &\in \lambda \partial g(x^{n_k+1}) + (x^{n_k+1} - (x^{n_k} - \lambda \stochgrad{i_k}{F}{x^{n_k}}{y^{n_k}} ))  \Leftrightarrow \\
    \lambda \grad F(x^{n_k+1}) - \lambda \grad F(x^{n_k+1}) &\in \lambda \partial g(x^{n_k+1}) + (x^{n_k+1} - (x^{n_k} - \lambda \stochgrad{i_k}{F}{x^{n_k}}{y^{n_k}} ))  \Leftrightarrow \\
    u^{n_k} &\in \partial g(x^{n_k+1}) + \grad F(x^{n_k+1})  \Leftrightarrow\\
    u^{n_k} &\in \partial (g + F)(x^{n_k+1}) \Leftrightarrow \\
    (x^{n_k+1}, u^{n_k}) &\in \gra{\partial (g + F)}
\end{align*}
{where $u^{n_k} = \lambda^{-1}(x^{n_k} - x^{{n_k}+1}) +  \grad F(x^{{n_k}+1}) - \stochgrad{i_k}{F}{x^{n_k}}{y^{n_k}}$. As ${n_k} \rightarrow \infty$, $x^{n_k} - x^{{n_k}+1} \rightarrow 0$ and $\grad F(x^{{n_k}+1}) - \stochgrad{i_k}{F}{x^{n_k}}{y^{n_k}} \rightarrow 0$. Thus, $u^{n_k} \rightarrow 0$ almost surely. In the second to last equivalence above, since $\partial F$ has full domain, we used the identity $\partial(g+F) = \partial g + \partial F$ by \cite[Corollary 16.48]{bauschke2017convex}. Let us assume that the subsequence converges to $\Bar{x}$, that is $x^{n_k} \rightarrow \Bar{x}$. Now by \cite[Corollary 25.5]{bauschke2017convex} since $\partial F$ has full domain, $\partial(g + F)$ is maximally monotone. Using \cite[Proposition 20.37.ii]{bauschke2017convex}, we get $(\Bar{x}, 0) \in \gra{\partial(g+F)}$ which implies that $0 \in \partial g(\Bar{x}) + \grad F(\Bar{x})$. This clearly means that all  sequential cluster points of $\seq{x^k}$ belong to $\mathcal{X}^\star$.} Now, invoking  \cite[Proposition 2.3.iv]{combettes2015stochastic}, implies that $\seq{x^k}$ converges almost surely to a $\mathcal{X}^\star$-valued random variable. Invoking Proposition~\ref{prop:p_pd_equivalence} concludes the proof.
\qed

\subsection{Proof of Theorem~\ref{thm:main}}\label{app:theorem5.2}

In the following proof, all the convergences and limits hold almost surely, even if it is not explicitly mentioned.

Let $I_{\mathrm{bm}}$ and $I_{\mathrm{aa}}$ be the sets of indices for which the next iterate is obtained by a \text{\basicmethod} step and an \text{\accalg} step, respectively. These index sets satisfy $I_{\mathrm{bm}}\cap I_{\mathrm{aa}}=\emptyset$ and $I_{\mathrm{bm}}\cup I_{\mathrm{aa}}=\mathbb{N}$. Note that if the cardinality of $I_{\mathrm{aa}}$ is finite, after a finite number of steps, the algorithm will reduce to the basic method, which we know is convergent by Theorem~\ref{thm:basicMethod}. Therefore, we assume that $|I_{\mathrm{aa}}|$ is infinite.

From Theorem~\ref{thm:basicMethod}, for all $k \in I_{\mathrm{bm}}$ and all $\z^\star \in \zer \Rbar$ we have
\begin{equation} \label{eq:SVRGseq}
    \Exp_k \left\Vert z^{k+1} - \z^\star \right\Vert_{\Gamma} ^2  \leq  \left\Vert z^k - \z^\star \right\Vert_\Gamma ^2 - \zeta_k,
\end{equation}
where $\zeta_k\geq 0$. Using the identity $\Exp \| X - \Exp X\|^2=\Exp \|X \|^2 - \|\Exp X \|^2$, we have $\|\Exp X \|^2\leq\Exp \|X \|^2$. Thus, from \eqref{eq:SVRGseq} and $\zeta_k\geq 0$ for all $k\in I_{\mathrm{bm}}$ we have
\begin{align*}
    \left(\Exp_k\norm{z^{k+1}-\z^\star}_{\Gamma}\right)^2 \leq \Exp_k\norm{z^{k+1}-\z^\star}_{\Gamma}^2  \leq  \norm{z^k-\z^\star}_\Gamma^2 - \zeta_k \leq \norm{z^k-\z^\star}_\Gamma^2.
\end{align*}
Therefore, for all $k\in I_{\mathrm{bm}}$ we have 
\begin{align}\label{eq:Exp_I_bm}
    \Exp_k\norm{z^{k+1}-\z^\star}_{\Gamma} \leq \norm{z^k-\z^\star}_\Gamma.
\end{align}
On the other hand for all $k\in I_{\mathrm{aa}}$, by the triangle inequality and for all $\z^\star \in \zer \Rbar$, we have
\begin{equation*}
    \norm{z^{k+1}-\z^\star}_\Gamma  \leq  \norm{z^k-\z^\star}_\Gamma + \norm{z^{k+1}-z^k}_\Gamma.
\end{equation*}
Using the safeguard condition~\eqref{eq:safeguard2} and that $\z^+=z^{k+1}$ for all $k\in I_{\mathrm{aa}}$, we obtain
\begin{equation} \label{eq:alpha_l}
   \norm{z^{k+1} - \z^\star}_\Gamma  \leq  \norm{z^k - \z^\star}_\Gamma + DV(z^{k}).
\end{equation}
Using the fact that  $\Exp_k\norm{z^{k+1}-\z^\star}_\Gamma=\norm{z^{k+1}-\z^\star}_\Gamma$ holds for all $k\in I_{\mathrm{aa}}$ since the acceleration method is deterministic, by combining \eqref{eq:Exp_I_bm} and \eqref{eq:alpha_l}, we conclude that 
\begin{equation} \label{eq:alpha_l_2}
    \Exp_k\norm{z^{k+1} - \z^\star}_\Gamma  \leq  \norm{z^k - \z^\star}_\Gamma + \sigma_k
\end{equation}
holds for all $k\in\nat$, where 
\begin{equation*}
    \sigma_k =
        \begin{cases}
            0     & k \in I_{\mathrm{bm}} \\
            D V(z^{k})  & k \in I_{\mathrm{aa}}
        \end{cases}.
\end{equation*}
Due to \eqref{eq:safeguard1}, $\seq{\sigma_k}$ is summable and $\seq{\norm{z^k-\z^\star}}$ converges a.s. \cite[Lemma 2.2]{combettes2015stochastic} and is therefore a.s. bounded. Next, by squaring both sides of \eqref{eq:alpha_l}, for all $k\in I_{\mathrm{aa}}$, we get
\begin{equation*} 
    \norm{\z^{k+1}-\z^\star}_\Gamma^2  \leq  \norm{\z^k-\z^\star}_\Gamma^2 + 2\norm{\z^k-\z^\star}_\Gamma  DV(z^{k}) + (DV(z^{k}))^2.
\end{equation*}
Defining $\beta_k := 2\norm{\z^k-\z^\star}_\Gamma DV(z^{k})+(DV(z^{k}))^2$ and using $\Exp_k\norm{z^{k+1}-\z^\star}_\Gamma=\norm{z^{k+1}-\z^\star}_\Gamma$ for all $k\in I_{\mathrm{aa}}$, we get for all $k\in I_{\mathrm{aa}}$ that
\begin{equation} \label{eq:alpha-squared}
    \Exp_k\norm{\z^{k+1}-\z^\star}_\Gamma^2  \leq  \norm{\z^k-\z^\star}_\Gamma^2 + \beta_k.
\end{equation}
Since we have concluded that $(\|\z^k - \z^\star\|_\Gamma)_{k\in I_{\mathrm{aa}}}$ is bounded a.s. and $(V(z^{k}))_{k\in I_{\mathrm{aa}}}$ is absolutely summable, $(\beta_k)_{k\in I_{\mathrm{aa}}}$ is a.s. absolutely summable as well. 
Combining \eqref{eq:SVRGseq} and \eqref{eq:alpha-squared} implies that
\begin{equation} \label{eq:OverallSeq1}
    \Exp_k\lv z^{k+1} - \z^\star\rv_\Gamma^2  + \nu_k \leq  \lv\z^k - \z^\star\rv_\Gamma^2  + \eta_k,
\end{equation}
where
\begin{equation*}
    \eta_k =
        \begin{cases}
            0     & k \in I_{\mathrm{bm}} \\
            \beta_k  & k \in I_{\mathrm{aa}}
        \end{cases},
\qquad{\hbox{and}}\qquad    \nu_k =
        \begin{cases}
            \zeta_k     & k \in I_{\mathrm{bm}} \\
            0   & k \in I_{\mathrm{aa}} 
        \end{cases}.
\end{equation*}
Therefore, by \cite[Proposition 2.3.i]{combettes2015stochastic}, $\seq{\nu_k}$ is summable. Now,  in \cite[Proposition 2.3]{combettes2015stochastic} setting $\phi:\z \mapsto \z^2$,  \cite[Proposition 2.3.iii]{combettes2015stochastic} implies that $\seq{\lv\z^k - \z^\star\rv^2_\Gamma}$  and evidently $\seq{\lv \z^k - \z^\star \rv}$ are convergent.

For the last part of the proof, fix $\z^\star \in \zer \Rbar$ and denote the set of sequential cluster points of $\seq{\z^k}$ by $\mathcal{C}$. Since $\seq{\|\z^k-\z^\star\|}$ is convergent, the sequence $\seq{\z^k}$  is bounded, and therefore, it has at least one convergent subsequence by the Bolzano--Weierstrass theorem. Denote this subsequence by $\seq{\z^{n_k}}$ and its associated sequential cluster point by $\z_c^\star=(x_c^\star,y_1^\star,\ldots,y_N^\star)$. As the problem is finite-dimensional, using Lipschitz continuity of the operator $\Rbar$ (Proposition~\ref{prop:Lipschitz-cont}) we have 
\begin{equation*}
    \norm{\z^{n_k}-\z_c^\star} \rightarrow 0 \qquad \Longrightarrow \qquad  \norm{\Rbar\z^{n_k}-\Rbar\z_c^\star} \rightarrow 0
\end{equation*}
which means that $ \Rbar \z^{n_k} \rightarrow \Rbar \z_c^\star $. Note that $(z^{n_k})_{k\in\nat}$ is constructed by the points that are generated by either the basic method or the acceleration algorithm. For the subsequence of points in $(z^{n_k})_{k\in\nat}$ that are obtained from the basic method, that is $(\z^{n_k+1})_{k\in{I_{\mathrm{bm}}}}$, since $\seq{\nu_{n_k}}$ is summable, so is $\seq{\zeta_{n_k}}$. Then, using the same approach as in the last part of the proof of Theroem~\ref{thm:basicMethod}, we can show that $(\Rbar z^{n_k+1})_{k\in{I_{\mathrm{bm}}}}$ converges to zero.  For the subsequence of the points in $(z^{n_k})_{k\in\nat}$ which are generated by the acceleration algorithm, that is $(\z^{n_k+1})_{k\in I_{\mathrm{aa}}}$, it is evident from the definition of the merit function in \eqref{eq:merit_func}, that convergence of $(V(\z^{n_k+1}))_{k\in{I_{\mathrm{aa}}}}$  to zero---which is dictated by the safeguard condition---enforces convergence of $(\Rbar \z^{n_k+1})_{k\in{I_{\mathrm{aa}}}}$ to zero as well. Therefore, for  $(\z^{n_k})_{k\in\nat}$ as a whole, we have $z^{n_k}\to0$ as $k\to\infty$. Then, it follows from $ \Rbar \z^{n_k} \rightarrow \Rbar \z_c^\star $ that $\Rbar \z_c^\star  = 0$. Thus, $\z_c^\star$ belongs to $\zer\Rbar$. The same implication can be made for all other sequential cluster points of $\seq{\z^k}$ which means that all sequential cluster points of $\seq{\z^k}$ belong to $\zer\Rbar$, that is  $\mathcal{C}\subset\zer\Rbar$. Finally, by  \cite[Proposition 2.3.iv]{combettes2015stochastic}, the sequence $\seq{\z^k}$ converges a.s. to a point $\Bar{z} \in \zer \Rbar$. Now, by Proposition~\ref{prop:dual-convergence}, $x^k \rightarrow \Bar{x}$ and for all $i$, $y_i^k \rightarrow \grad f_i (\Bar{x})$ a.s., where $\Bar{x}$ is the solution of problem~\eqref{eq:optProblem-repeat}. By this, the proof is complete.
\qed

\section*{Appendix B}
\renewcommand{\thesubsection}{B.\arabic{subsection}}
\setcounter{subsection}{0}
\renewcommand{\thealgorithm}{B.\arabic{algorithm}}
\setcounter{algorithm}{0}

\subsection{Anderson acceleration} \label{subsec:Anderson}
Anderson acceleration can be exploited to accelerate convergence of the the fixed--point iteration of the form
\begin{equation*}
    x^{k+1} = \T(x^{k})
\end{equation*}
where $\T:\reals^n\rightarrow\reals^n$ is either a contraction or an averaged operator. A variant of Anderson acceleration,  which is equipped with Tikhonov regularization on its inner least--squares problem, is given in Algorithm~\ref{alg:AA} \cite{scieur2016regularized}.

\begin{algorithm}
	\caption{Anderson Acceleration}
	\begin{algorithmic}[1]
	    \State \textbf{input:} $y^0\in\reals^{d}$, $m\geq1$.
		\For {$k=0,1,2,\ldots$}
		    \State set $m_k = \min\{m,k\}$
		    \State find the iterate using the fixed-point map $ x^{k} = \T(y^{k})$
		    \State form $\mathscr{R}^k = (r^{k-m_k},\ldots,r^k)$ where $r^{j} = y^{j} - x^{j}$ for $j\in\{k-m_k,\ldots,k\}$
		    \State determine $\alpha^{(k)} = (\alpha_0^{(k)},\ldots,\alpha_{m_k}^{(k)})$ that solves
		        \begin{equation*}
                \begin{aligned}
                   &\underset{~~\alpha^{(k)}\in\reals^{m_k+1}}{\mathrm{minimize}} ~~~ \norm{\mathscr{R}_n\alpha^{(k)}}_2^2 + \xi_k\norm{\alpha^{(k)}}_2^2 \\
            		            &\mathrm{~~subject ~ to ~~~~}  \mathbf{1}^T\alpha^{(k)}=1
                \end{aligned}
                \end{equation*}
		    \State  $y^{k+1}=\sum_{i=0}^{m_k}\alpha_{i}^{(k)}x^{k-m_k+i}$
		\EndFor
	\end{algorithmic}
\label{alg:AA}
\end{algorithm}

\subsection{Limited--memory BFGS}\label{subsec:lBFGS}
If the objective function of a convex optimization problem is twice continuously differentiable, an effective way of solving it, is to use quasi--Newton methods. One of the most well-known quasi--Newton methods is the limited--memory BFGS (lBFGS) which has been vastly used in many areas. lBFGS is a variant of BFGS method that uses a limited amount of computer's memory and in that sense is cheaper than its parent, BFGS method. Hence, unlike BFGS algorithm, its limited--memory version can be used to solve large--scale problems. The lBFGS method can be stated as in Algorithm \ref{alg:lBFGS} \cite{nocedal2006numerical}.

\begin{algorithm}
	\caption{limited--memory BFGS}
	\begin{algorithmic}[1]
	    \State  Define: $s^k = x^{k+1} - x^k$, $u^k = \grad f(x^{k+1}) - \grad f(x^k)$ and $\rho_k = ((u^k)^T s^k)^{-1}$.
        \State \textbf{input:} $x^0$ and the memory stack size $m \geq 1$.
		\For {$k=1,2,\ldots$}
		    \State $H_0^k = \tfrac{(s^{k-1})^T u^{k-1}}{(u^{k-1})^T u^{k-1}}I$
		    \State $q = \grad f(x^k)$
		    \For {$i = k-1,\hdots, \text{min}\{k-m,0\}$}
		        \State $\alpha_i = \rho_i (s^i)^T q$
		        \State $q = q - \alpha_i u^i$
		    \EndFor
		    \State $r = H_0^k q$
		    \For {$i = \text{min}\{k-m,0\},\hdots,k-1$}
		        \State $\beta = \rho_i (u^i)^T r$
		        \State $r = r + s_i(\alpha_i - \beta)$
		    \EndFor
		    \State $p^k = r$
		    \State \text{compute }$x^{k+1} = x^k - \lambda_k p^k$, \text{ where $\lambda_k$ is to satisfy a line search condition}
		    \If{$k>m$}
		        \State \text{discard $s^{k-m}$ and $u^{k-m}$}
		        \State \text{compute and save $s^{k}$ and $u^{k}$}
		    \EndIf
		\EndFor
	\end{algorithmic}
\label{alg:lBFGS}
\end{algorithm}

\printbibliography


\end{document}